\documentclass[10pt]{amsart}
\usepackage[margin=1in]{geometry}
\usepackage{amsmath,amssymb,mathtools}
\usepackage[dvipsnames]{xcolor}
\input{eo_cplx_VBS.sty}
\makeatletter
\renewcommand{\email}[2][]{%
  \ifx\emails\@empty\relax\else{\g@addto@macro\emails{,\space}}\fi%
  \@ifnotempty{#1}{\g@addto@macro\emails{\textrm{(#1)}\space}}%
  \g@addto@macro\emails{#2}%
}
\makeatother

\title[Enumerating stably trivial vector bundles]
{Enumerating stably trivial vector bundles \\ with higher real $K$-theory}
\author{Hood Chatham}
\author{Yang Hu}
\author{Morgan Opie}
\address[H. Chatham]{Cloudflare}
\email{roberthoodchatham@gmail.com}
\address[Y. Hu]{University of Regina}
\email{yang.hu@regina.ca}
\address[M. Opie]{Northwestern University}
\email{mpopie@northwestern.edu}

\date{}              

\begin{document}

\begin{abstract}
This paper explores periodic phenomena in the group $\operatorname{Vect}_r^0(\mathbb{CP}^{r+c})$ of stably trivial, complex rank $r$ topological vector bundles on $\mathbb{CP}^{r+c}$. For $1 \leq c < r$ and $c\leq 2p-3$, we give a complete computation of the $p$-torsion in $\operatorname{Vect}_r^0(\mathbb{CP}^{r+c})$, and we relate these $p$-torsion bundles to vector bundles on spheres. We also compute $\operatorname{Vect}_r^0(\mathbb{CP}^{r+c})$ in full when $c=3$, extending computations of the second-named author when $c=1$ and $c=2$. Finally, for a fixed corank $c$ which is larger relative to the prime $p$, we show there are families of $p$-torsion in $\operatorname{Vect}_r^0(\mathbb{CP}^{r+c})$ that are periodic in $r$ with period $p$. We detect these using Weiss calculus and certain higher real $K$-theory homology groups.
\end{abstract}
\maketitle
\tableofcontents

\section{Introduction}

The study of complex topological vector bundles on compact manifolds has a rich history. Some of the most-studied examples include vector bundles on spheres and vector bundles on complex projective spaces. 

Complex vector bundles on spheres are relatively well-understood in the {\em metastable range}, due to work of Bott, Kervaire, Matsunaga, Mimura, Toda, and others; these results are succinctly summarized in \cite[Section 3.2]{Mimura_HBAT}. Many of these results build on \cite[Theorem 4.3]{Toda59}, which states that, for $\frac{n}{2} \leq r <n$, there is an isomorphism \begin{equation}\label[iso]{iso:unitarymetastable}\pi_{2n}BU(r) \cong \pi_{2n}^{\op{st}} \Sigma \cp{r}{n}\end{equation} where $\cp{r}{n}$ is the space obtained by collapsing cells in $\cp{}{n}$ of dimension less than $2r$ to a point, and $\pi_*^{\op{st}}$ denotes stable homotopy. 
This allows for computation of unstable homotopy of unitary groups in terms of the {\em stable} (and hence more computable) homotopy groups of well-understood cell complexes. Moreover, it highlights patterns in unstable homotopy groups of unitary groups for fixed corank $c=n-r$, provided $\frac{n}{2} \leq r <n$: complex corank $c$ vector bundles on $S^{2n}$ are periodic in $n$, with period the Atiyah--Todd number $M_{n-r+1}$.\footnote{This follows from observing that $\cp{r}{n}$ is the Thom space of the sum of $r$ copies of the tautological line bundle on $\cp{}{n-r}$ and applying James periodicity; see \Cref{rm:James} for more discussion.}

In contrast, until recently, complex topological vector bundles on complex projective spaces were studied in an ad hoc manner. This changed with work of the second-named author, which gives an upgrade of \Cref{iso:unitarymetastable} in terms of Weiss' {\em unitary calculus} \cite[Proposition 2.3]{Hu}. An explicit consequence of this is the following: let $\Vect_r^0(\cp{}{n})$ denote the set of
isomorphism classes of stably trivial rank $r$ bundles on $\cp{}{n}$. For  $\frac{n}{2}\leq r \leq n-1$, there is a natural bijection 
\begin{equation}\label[iso]{eq:Yang-thm}\Vect^0_r(\cp{}{n})\cong \{\cp{r}{n},\Sigma \cp{r}{n}\},\end{equation} where $\{-,-\}$ denotes stable homotopy classes of maps \cite[Theorem 2.1]{Hu}. 

\begin{conv} For the rest of the introduction, we always assume that $n$ and $r$ are metastable, meaning that $\frac{n}{2}\leq r < n$, so that \Cref{eq:Yang-thm} holds. \end{conv}

This set of stable maps on the right-hand side of \Cref{eq:Yang-thm} exhibits many nice properties. For fixed corank $c=n-r$, James periodicity implies that the enumeration of stably trivial rank $r$ bundles on $\cp{}{n}$ is periodic in $r$, with period the Atiyah--Todd number $M_{c+1}$
\cite[Theorem 1.2]{AdamsWalker}. 
\begin{rmk}\label{rm:James} A formula for $M_k$ is given in \cite[Page 100]{AdamsWalker}. Small values are $M_2=2$, $M_3=M_4=24$, $M_5=2880$. $M_{k}$ grows rapidly as $k$ grows -- for instance,
\[\nu_p(M_k) \geq \left\lfloor \frac{k-1}{p-1} \right\rfloor\]
gives a conservative bound for $p$-divisibility of $M_k$. Here, $\nu_p$ denotes $p$-adic valuation. In particular, the periodicity grows rapidly with $c$.
\end{rmk}
In \cite[Theorem 1.1 and 1.2]{Hu}, the second-named author uses \Cref{eq:Yang-thm} together with the observed periodicity and Adams spectral sequence computations to
enumerate all stably trivial, non-trivial metastable complex bundles on complex projective spaces in corank $c=1$ and corank $c=2$.  One result in this paper is to complete the corank $c=3$ case:

\begin{thmx}[{See \Cref{thm:corank3-3} and \Cref{thm:corank3-2}}] \label{thm:corank3-intro}
Let $r \geq 3$. The number of isomorphism classes of stably trivial rank $r$ bundles on $\cp{}{r+3}$ is as follows:

\[
\begin{tikzpicture}
\matrix (m) [matrix of math nodes,
nodes in empty cells,nodes={minimum width=2ex,
minimum height=2ex,outer sep=-1pt},
column sep=.2ex,row sep=.1ex]{
r \mod{24} & 0 & 1 & 2 & 3 & 4 & 5 & 6 & 7 & 8 & 9 & 10 & 11 & \, \\
\#\Vect^0_r(\cp{}{r+3}) & 24 & 4 & 3 & 12 & 8 & 12 & 3 & 8 & 24 & 12 & 1 & 12 & \, \\
\\
r \mod{24} & 12 & 13 & 14 & 15 & 16 & 17 & 18 & 19 & 20 & 21 & 22 & 23 & \, \\
\#\Vect^0_r(\cp{}{r+3}) & 24 & 4 & 3 & 24 & 8 & 12 & 3 & 4 & 24 & 12 & 1 & 24 \strut \\
};
\draw[thick] (m-1-1.north west) -- (m-1-13.north east);
\draw[] (m-1-1.south west) -- (m-1-13.south east);
\draw[thick] (m-4-1.north west) -- (m-4-13.north east);
\draw[] (m-4-1.south west) -- (m-4-13.south east);
\end{tikzpicture}
\]

\end{thmx}

\begin{rmk} Theorem A is obtained by combining the $2$- and $3$-primary calculations. There is no $p$-primary contribution for $p\geq 5$, an observation that is justified in \Cref{thm:B} below. \end{rmk} 

The optimistic reader might hope to compute the group $\Vect_r^0(\cp{}{n})$ -- or at least its order -- in full. However, at least by hand, this seems implausible: for example, for $c= 4$, there are $M_5=2^6\cdot 3^2\cdot 5$ cases to consider. Rather than pursuing a complete computation, we turn our attention to certain patterns.

An abstract consequence of \Cref{eq:Yang-thm} is that stably trivial rank $r$ vector bundles on $\cp{}{n}$ actually form a finite abelian group (finiteness can be deduced from a classical obstruction theory argument). Of course, such a group is a direct sum of its $p$-power torsion subgroups, and so any stably trivial vector bundle can be decomposed into $p$-primary constituents. Understanding this group structure as an explicit construction on vector bundles would be ideal, but at present it remains mysterious. However, a homotopy-theoretic compromise is to describe the $p$-power torsion one prime at a time, and return to the problem of understanding the group structure later. 

For fixed corank $c=n-r$, the primes can be grouped roughly into four regimes, relative to the corank:
\begin{itemize}
\item[i.] {\bf Large primes}: $c<p-1$. These do not contribute to $\Vect_r^0(\cp{}{n})$ (see \Cref{thm:B} below).
\item[ii.] {\bf Intermediate primes}: $p-1 \leq c \leq 2p-3$. These contribute only $p$-torsion to $\Vect_r^0(\cp{}{n})$, in a way that we understand completely (see \Cref{thm:B}), and that we can analyze geometrically (see \Cref{thm:Y}). 
\item[iii.] {\bf Small primes}: $2p-3<c<2p^2-p-2$. We have a limited understanding of $p$-power torsion for these primes, but can detect infinite families of partial information (see \Cref{thm:A} and \Cref{thm:X}) and explore patterns via explicit computation for specific values of $p$ and $c$ (e.g., \Cref{thm:corank3-intro} for $p=2$ and $c=3$, and \Cref{thm:corank_2p-2} below for $p=3$ and $c=4$).
\item[iv.] {\bf Very small primes}: $c\geq 2p^2-p-2$. We generally do not discuss these primes in this paper.
\end{itemize}
We describe this behavior more precisely in the next few results, beginning with cases (i) and (ii).
\begin{thmx}[See \Cref{cor:bundle_ct1}]\label{thm:B} Let $p$ be a prime. Let $r$ be a given rank and $n$ a given dimension, and suppose that $\frac{n}{2}\leq r<n$, so that rank $r$ bundles on $\cp{}{n}$ are metastable. Let $c=n-r$.
\begin{itemize}
\item  If $c<p-1$, then $\Vect^0_r(\cp{}{n})$ contains no $p$-power torsion.

\item If $p-1\leq c \leq 2p-3$, the $p$-power torsion in $\Vect_r^0(\cp{}{n})$ is isomorphic to $\Z/p$, if $r+i\equiv 0 \pmod p$ for some $i\in \{ 0,\ldots,c-(p-1)\}$, and is zero otherwise.
\end{itemize}
\end{thmx}
\Cref{thm:B} gives a small number of $p$-torsion vector bundles for each $p$ with $p-1\leq c\leq 2p-3,$ so we might hope to explicitly describe such bundles in a geometric way. We show that essentially all $p$-torsion bundles in this range are pulled back from spheres:
\begin{thmx}[See \Cref{geometric-easy}]\label{thm:Y} Let $p$ be prime and let $n,r$ be positive integers. Suppose that $n-r\leq r <n$. Let $c=n-r$ be such that $p-1\leq c \leq 2p-3$. Then non-trivial, stably trivial $p$-torsion rank $r$ vector bundles on $\cp{}{n}$ are of two types:
\begin{itemize}
\item[]
\begin{itemize}
\item[{\bf Type 1.}] If $n\equiv -1 \pmod{p}$ and $\Vect_r^0(\cp{}{n})$ contains a non-zero $p$-torsion bundle $V$, then $V$ is pulled back from a non-trivial complex rank $r$ vector bundle on $S^{2n}$ under projection onto the top cell $\cp{}{n} \to S^{2n}$. 
\item[{\bf Type 2.}] If $n\not \equiv -1 \pmod{p}$ and $\Vect_r^0(\cp{}{n})$ contains a non-zero $p$-torsion bundle $W$, then $W$ is the unique extension of a bundle of Type 1 over a cellular inclusion $\cp{}{k}\to \cp{}{n}$ for $k\equiv -1 \pmod{p}$.
\end{itemize}
\end{itemize}
\end{thmx} 

From these computations, one can also observe that all non-trivial, stably trivial $p$-power torsion bundles in the range $p-1\leq c \leq 2p-3$ are, in a precise sense, detected by a generalized cohomology theory called {\em higher real $K$-theory}, which we write as $\eo_{p-1}$ given the implicit prime $p$. The relevance of $\eo_{p-1}$ is not entirely surprising: the story of classifying unstable vector bundles on complex projective spaces begins with Atiyah and Rees' seminal work on complex rank $2$ topological vector bundles on 
$\cp{}{3}$ \cite[Theorem 2.8]{AR}. In their paper, Atiyah--Rees use real $K$-theory to distinguish non-isomorphic rank $2$ vector bundles with the same $K$-theory class. Real $K$-theory is closely related to the higher real $K$-theory $\eo_{1}$ at the prime $2$. Further evidence of the relevance of higher real $K$-theories comes from the third author's work classifying rank $3$ bundles on $\cp{}{5}$ \cite[Theorem 1.4]{MO} using the generalized cohomology theory of topological modular forms localized at the prime $3$. The latter theory is closely related to higher real $K$-theory $\eo_2$ at the prime $3$. So, viewed together, \cite{AR} and \cite{MO} suggest that the simplest topological vector bundles that {\em cannot} be distinguished by complex topological $K$-theory {\em can} be teased apart using higher real $K$-theories.

Our next objective is to run with this idea: exactly how much of the group $\Vect_r^0(\cp{}{n})$ can we see through the lens of higher real $K$-theories? To address this question, we prove an abstract detection theorem. 
\begin{thmx}[See \Cref{prop:surjective}]\label{thm:A}
Let $p$ be an odd prime and $r,n$ satisfy $\frac{n}{2}\leq r <n$. Let $\eop$ denote the generalized homology theory of higher real $K$-theory of height $p-1$ at the implicit prime $p$. For $n-r< 2p^2-p-2$, there is a surjection \begin{equation}\label[group]{eq:eo-homology}\Vect_r^0(\cp{}{n})\to (\eop)_{0}(Z_{n,r}),\end{equation}
where $(\eop)_{0}(Z_{n,r})$ is the zeroth generalized homology group of a certain well-understood finite spectrum with respect to $\eop$.
\end{thmx}
\begin{rmk} Explicitly $Z_{n,r} \cong \Sigma \cp{r}{n} \otimes D\cp{r}{n}$, where $\cp{r}{n}$ is a stunted projective spectrum, $\otimes$ denotes the smash product of spectra, and $D\cp{r}{n}$ denotes the Spanier--Whitehead dual of $\cp{r}{n}$.\end{rmk}

\begin{rmk} 
For fixed corank $c=n-r$, the group $\eop_0(Z_{n,r})$ is periodic in $n$, with period $p$. This contrasts with $\Vect_r^0(\cp{}{n})$, whose period is typically much larger and depends strongly on $c$ (see \Cref{rm:James}). \end{rmk}

As a consequence of \Cref{thm:A}, we can show that non-trivial, stably trivial vector bundles exist by producing non-zero elements in the mysterious but computationally tractable group $(\eop)_{0}(Z_{n,r})$. 

To complement this result, we give a complete computation of the group $(\eop)_{0}(Z_{n,r})$ in the range $c<\frac{2p^2-p-2}{2}$. This result is perhaps of greatest interest to those with a background in computational chromatic homotopy theory, and is technical to state, so we refer the reader to \Cref{thm:total} below for the full result. 
However, we note a concrete consequence here:
\begin{thmx}[See \Cref{cor:bigcount}] \label{thm:X} Let $p$ be an odd prime and $i\geq 0$.
The number of isomorphism classes of stably trivial rank $p^2+ip$ bundles on $\cp{}{(p-1)^2+p^2+ip}$ is a nonzero multiple of  $p^{ p-2}$.
\end{thmx}

The existence of nontrivial families of bundles as in \Cref{thm:X} raises many questions. It would be of great interest to give a geometric construction of such bundles, along the lines of \Cref{thm:Y}. We give some partial results and some negative results along these lines, by showing that many of the vector bundles that we detect are {\em not} pulled back from spheres (see \Cref{lem:same-adams}). Thus, explicitly constructing these vector bundles requires some more complicated geometry.

\begin{rmk} 
We conclude with a discussion of periodic phenomena and how they should be interpreted. The periodicity implied by \Cref{eq:Yang-thm} is genuinely present in the structure of stably trivial vector bundles, although we do not have an entirely satisfactory geometric explanation for it. In the regime $p-1\leq n-r\leq 2p-3$, \Cref{thm:Y} gives a reasonable way to relate these periodic vector bundles on complex projective spaces to periodic phenomena for vector bundles on spheres. On the other hand, the $p$-periodic families of $p$-torsion vector bundles detected in higher real $K$-theories using \Cref{thm:A} (e.g., as in \Cref{thm:X}) are algebraic in nature. It seems reasonable to hope that the patterns detected in higher real $K$-theory may be avatars of geometric structure to be determined. \end{rmk} 


\subsection{Methods}\label{methods} 
We now give a more technical overview of the main tools used in our paper. 

\Cref{eq:Yang-thm} converts a question about topology into a more tractable question about stable homotopy
groups. The first results in this paper are a direct calculation of the $p$-primary part of the group of stable maps
$\{ \cp{r}{n},\Sigma\cp{r}{n}\}.$ To establish an analogy, we note that this group of stable maps might also be written as
$ \hmapsp{\cp{r}{n}}{\Sigma\cp{r}{n}}$, i.e., homotopy classes of maps of $\sphere$-modules.

For arbitrary $n$ and $r$, the Adams spectral sequence computing $p$-primary stable maps from $\cp{r}{n}$ to $\Sigma \cp{r}{n}$ is quite complicated, so we might simplify our approach by base-changing from the category of $\sphere$-modules (i.e., the category of spectra) 
 to the category of $R$-modules where $R$ is chosen to simplify the
problem. We will usually base change to $R$-modules for $R$ a higher real $K$-theory spectrum of height $p-1$ at an
implicit prime $p$.  However, passing from $\sphere$-modules to $\eop$-modules need not generally preserve or reflect information about $\sphere$-module maps. As such, the detection result \Cref{thm:A} is necessary. We can rephrase \Cref{thm:A} as stating that a Hurewicz map 
\begin{equation}\label[empty]{Hurw-1}\hmapsp{\cp{r}{n}}{\Sigma\cp{r}{n}} \to \hmaps{\eop}{\eop \otimes \cp{r}{n}}{ \eop \otimes \Sigma \cp{r}{n}} \end{equation}
is surjective for a certain range of $r$ and $n$.

\begin{rmk}
The higher real $K$-theory spectra $\eop$ are the
simplest in a hierarchy of {\em chromatic theories} built to detect certain constituents of the stable homotopy groups of spheres. 
For each prime $p$, and integer $m>1$, there is a corresponding height $m(p-1)$ higher real $K$-theory spectrum $\eo_{m(p-1)}$. Studying the $\eo_{m(p-1)}$-homology of $\cp{r}{r+c} \otimes D\cp{r}{r+c}$ as $r$, $c$, and $m$ vary seems likely to give extensive information about stably trivial vector bundles on complex projective spaces. However, computations with the spectra $\eo_{m(p-1)}$ for $m>1$ are extremely challenging. Thus, $\eop$-theory is a Goldilocks spectrum for our problem: relevant $\eop$-homology is highly computable, while also detecting interesting structure. \end{rmk}

To prove \Cref{thm:A}, we appeal to the theory of algebraic $\eop$-modules developed
by the first author \cite{Chat}. A spectrum $Z$ has algebraic $\eop$-theory if
$\eop\sm Z$ splits as a sum of simpler $\eop$-modules in a manner that lifts the
decomposition of the cohomology groups $H^*(Z;\mathbb F_p)$ into indecomposable $P^1$-modules (here, $P^1$ denotes the $p$-primary Steenrod power operation of degree $2(p-1)$). Many
naturally occurring spectra have algebraic $\eop$-theory, including complex
projective spaces, stunted complex projective spaces, and their smash products.

The proof of  \Cref{thm:A} can be broken into three main steps.
First, we construct a cofiltration of an arbitrary even-celled finite spectrum $Z$ by particularly simple spectra which we call $X_l$, which are defined for each $l\in \{1,\ldots, p\}$ and characterized by their mod $p$ cohomology (see \Cref{def:algeo}).

\begin{thmx}[{See \Cref{cor:exists_Xi_filtration}}]\label{prop:egegeg} Let $Z$ be an even, finite spectrum. Then $Z$ admits a cofiltration 
\begin{equation}\label[diagram]{eq:filtfilt0}
\begin{tikzcd}[row sep=1.3em]
A_0 \arrow[d,"f_0"]& A_1  \arrow[d,"f_1"] & A_2 \arrow[d,"f_2"] & & A_n \ar[d,"\simeq"]
\\ 
Z=C_0\ar[r,"g_0"]& C_1 \ar[r,"g_1"]& C_2 \ar[r,"g_2"]& \cdots \,\,  \ar[r,"g_{n-1}"]&C_n&,
\end{tikzcd}\end{equation}
where $A_i$ is the fiber of the $g_i$,
such that for all $i$ 
 \[A_i\simeq \Sigma^{s_i}X_{l_i}.\]  If additionally $Z$ has algebraic $\eop$-theory, then each $f_i\:A_i \to C_i$ splits after tensoring with $\eop$.
\end{thmx}
\begin{rmk} Since $\cp{r}{n} \otimes D \cp{r}{n}$ is an even-celled finite spectrum with algebraic $\eop$-theory, \Cref{prop:egegeg} gives us an $\eo$-split cofiltration of $\cp{r}{n} \otimes D \cp{r}{n}$.\end{rmk}

The second step in the proof of \Cref{thm:A} is to reduce the surjectivity of the Hurewicz map \Cref{Hurw-1}
to properties of the Hurewicz image for cofiltration quotients of the form $\Sigma^{2s} X_l$, appearing in \Cref{eq:filtfilt0} for $Z=\cp{r}{n} \otimes D\cp{r}{n}$. The final step is to compute the Hurewicz maps $\pi_* X_l \to (\eop)_*X_l$ in a suitable range.

To actually obtain lower bounds for vector bundles, we use tools from \cite{Chat} to compute $\eo_{-1}\cp{r}{n} \otimes D\cp{r}{n}$ in full provided $n-r<\frac{2p^2-p-2}{2}$. This is an involved combinatorial process.


\subsection{Paper outline}
\label{subsec:outline}

In \Cref{sec:background}, we review some history and survey various homotopical tools used in our enumeration results.
In \Cref{sec:first}, we give computations for $p$-completed corank $c$ bundles on
$\cp{}{r+c}$ for $1\leq c <r$ and $c$ small relative to $p$. We begin \Cref{subsec:small_corank1} by enumerating corank $p-1$, which is then generalized to $c\leq 2p-3$ in \Cref{subsec:small_corank2}. This proves \Cref{thm:B}.

 In \Cref{subsec:corank2p2}, we give a complete computation of corank $3$ vector bundles on $\cp{}{n}$, for $n$ at least $6$. We also give a computation of the $3$-primary part of $\Vect_r^0(\cp{}{r+4})$, illustrating that patterns become more complicated at corank $2p-2$.

We first explicitly use higher real $K$-theories in \Cref{sec:eo_detection}. In
\Cref{cofiltrations} we construct a cofiltration of
$\emod{\sphere}\left(\cp{r}{r+c},\Sigma \cp{r}{r+c}\right)$ by suspensions of
particularly nice spectra (see \Cref{prop:egegeg} above). In \Cref{computations_Xl}, we reduce
\Cref{thm:A} to a computation on the filtration quotients and
complete its proof. In \Cref{lowerbound_computation}, we apply \Cref{thm:A} to detect nontrivial rank $r$ bundles on
$\cp{}{r+c}$ when $c< \frac{2p^2-p-2}{2}$, with a complete computation of $\eo_{-1}(\cp{r}{r+c} \otimes D\cp{r}{r+c})$ (cf. \Cref{thm:total}). The consequences include \Cref{thm:X}.


\Cref{combinatorics} includes computations needed in
\Cref{lowerbound_computation} for the proof of \Cref{thm:total}. This section is combinatorial.

\subsection{Conventions}\label{subsec:conventions}
\begin{itemize}
\item Given pointed spaces $X$ and $Y$, we write $[X,Y]$ for homotopy classes of pointed maps from $X$ to $Y$.
\item Given spectra $X$ and $Y$, we write $\hmapsp{X}{Y}$ for homotopy classes of spectra maps from $X$ to $Y$. 
\item We will frequently consider Morava $E$-theory, and higher real $K$-theory spectra, of height $p-1$ at a prime $p$. These are usually abbreviated by $\e_{p\,\text{-}1}$ and $\eop$, respectively. When there is no potential for ambiguity, the spectrum $\eop$ will be abbreviated by $\eo$
 and $\e_{p\,\text{-}1}$ by $\e$.
\item Given $\eo$-modules $A$ and $B$, we write $\eomod(A,B)$ for the space of $\eo$-module maps from $A$ to $B$.  
\item We write $\otimes$ for the smash product of spectra.
\item Given a finite spectrum $X$, we write $DX$ for the Spanier--Whitehead dual of $X$.
\item We write $H^*$ for $H^*(-;\fp)$ and $H_*$ for $H_*(-;\fp)$.
\item In \Cref{sec:first} and beyond, an implicit prime $p$ is fixed. All spectra are implicitly completed at $p$, and all homotopy groups and groups of stable vector bundles are understood to be up to $p$-completion. Thus, in \Cref{sec:first} and beyond, $\hmapsp{X}{Y}$ denotes homotopy classes of maps of spectra from the $p$-completion of $X$ to the $p$-completion of $Y$. (See \Cref{subsec:background_loc} for discussion of $p$-completion.) 

\end{itemize}

\subsection{Acknowledgements}
The authors would like to thank Dev Sinha for his encouragement of this project, including travel funding provided under the Simons 
Foundation Collaboration Grant 422618. The authors also benefited greatly from
the opportunity to collaborate during South Central Topology Conference III.

The authors are also grateful to Prasit Bhattacharya, Mike Hill, Mike
Hopkins, Haynes Miller, Ross Staffeldt, and Alexander Smith for helpful discussions.

While working on this project, H.C. and M.O. were supported by the National Science Foundation under Award
Numbers 2002087 and 2202914, respectively.

\section{Background}
\label{sec:background}
\subsection{Weiss calculus and metastable vector bundles}\label{subsec:weiss}

Weiss calculus provides a framework to study continuous functors from the category of
finite-dimensional complex inner product spaces (and their linear isometric inclusions) to the
category of pointed topological spaces \cite{Weiss}. For each such functor $F$, Weiss calculus produces a tower of fibrations
\[\cdots \longrightarrow T_nF \longrightarrow T_{n-1}F 
\longrightarrow \cdots \longrightarrow T_1F \longrightarrow T_0F\]
known as the {\em Weiss tower} of $F$. For each
$n\geq 0$, $T_nF$ is an {\em $n$-polynomial} functor in the sense of Weiss \cite[Definition 5.1]{Weiss},
admitting a universal comparison from $F$ compatible with the tower.
The tower is said to {\em converge} to $F$ if $F \rightarrow \holim_n
T_nF$ is a pointwise weak equivalence.
The homotopy fiber $L_nF:= \text{hofib } (T_nF \rightarrow T_{n-1}F)$, known as the {\em $n$-th layer} of $F$, is an
{\em $n$-homogeneous} functor. A fundamental theorem of Weiss \cite[Theorem 7.3]{Weiss} classifies $n$-homogeneous 
functors as determined by (naive) $U(n)$-spectra. 
The bottom stage $T_0F$ of the Weiss tower is by definition 
\[(T_0F)(V) = \hocolim_k F(V\oplus \bb{C}^k),\] 
which is a constant functor and manifestly a stabilization.

Consider the functor $F=BU(-)$, given on objects by $V \mapsto BU(V)$. 
In this case, $T_0\left( BU(-) \right)$ is the constant 
functor $V\mapsto BU$ whose target classifies stable bundles. 
Higher layers of the tower have been studied in detail by 
Arone \cite[Theorem 2 and Theorem 4]{Arone02}. Arone's results imply that, if $r$ is the complex dimension of 
$V$ and $n\geq 2$, the space $L_nBU(V)$ is at least $4r$-connected.
This connectivity estimate can be applied to detect unstable complex vector bundles over a finite cell complex $X$. If the dimension of 
$X$ is no greater than $4r$, then $[X, BU(V)]$ is in natural bijection with $[X, T_1BU(V)]$. This is the so-called {\em metastable range}, 
and a rank $r$ bundle over $\cp{}{n}$ is metastable if $\frac{n}{2} \leq r < n$.  In the metastable range, stably trivial rank $r$ vector bundles on $\cp{}{n}$ admit the following identification (see \cite[Theorem 2.1]{Hu} for details):
\[\Vect_r^0(\cp{}{n}) \cong \hmapsp{ \cp{}{n}}{\Sigma \cp{r}{\infty}}\cong \hmapsp{ \cp{r}{n}}{ \Sigma \cp{r}{n}}.\]

\subsection{Stable homotopy theory completed at
\texorpdfstring{$p$}{p}}\label{subsec:background_loc} The fracture square in stable homotopy theory is analogous to the Hasse principle in number theory (namely checking an arithmetic property over $\bb{Z}$ is equivalent to checking it over $\bb{Q}$ and over the $p$-adic integers $\bb{Z}_p$ for every prime $p$). 
Recall that given a spectrum $Z$, its arithmetic completion or localization can be phrased as the Bousfield localization $L_E Z$ \cite{Bousfield} with respect to an appropriate Moore spectrum $E = MG$. Standard examples include:
\begin{itemize}
\item When $G = \bb{Q}$, $L_EZ$ is the rationalization of $Z$, written as $Z_{\bb{Q}}$.
\item When $G = \bb{Z}/p$ and when $Z$ is connective, $L_EZ$ is the $p$-completion of $Z$, denoted by $Z_p$.
\item When $G=\Z_{(p)}$, the integers localized at $p$, $L_EZ$ is the localization of $Z$ at the prime $p$.
\end{itemize}

\begin{thm}[{\cite[Proposition 2.9]{Bousfield}}]
   \label{thm:fracture_spectra}
   The following is a homotopy pullback square for every spectrum $Z$.
\begin{equation} \label[diagram]{fracture-square}
      \begin{tikzcd}[column sep = large]
           Z \ar[r] \ar[d] & \prod_p Z_p \ar[d] \\
           Z_{\bb{Q}} \ar[r] & (\prod_p Z_p)_{\bb{Q}}
      \end{tikzcd}
\end{equation}
\end{thm}
\Cref{fracture-square} is called the arithmetic square or fracture square, and allows us to work one prime at a time. 

\begin{lemma} 
   If $Z$ is a finite spectrum, then $\pi_*(Z_p)$ is the $p$-completion of $\pi_*(Z)$ at $p$.
\end{lemma}
This immediately implies the following.
\begin{cor}\label{cor:pcomplete}
Let $Y$ and $Z$ be finite spectra. Then the $p$-completion of $\hmapsp{Y}{Z}$
is isomorphic to $\hmapsp{Y_p}{Z_p}$. In particular, 
if $\hmapsp{Y}{Z}$ is a torsion group, $\hmapsp{Y_p}{Z_p}$ 
is identified with the $p$-power torsion in $\hmapsp{Y}{Z}$.
\end{cor}

Thus, to calculate $p$-divisibility of the group $\Vect_r^0(\cp{}{n})$ with $n,r$ in the metastable range,
we can compute homotopy classes of maps between the relevant $p$-complete spectra. 

\subsection{
    Higher real \texorpdfstring{$K$}{K}-theories of height
    \texorpdfstring{$p-1$}{p−1} at a prime
    \texorpdfstring{$p$}{p}
}\label{subsec:background_EO} 

Morava $E$-theories are fundamental in chromatic homotopy theory. 
Let $\bb{F}$ be a perfect field of characteristic $p$, and let
$\Gamma$ be a formal group of height $n$ over $\bb{F}$. The universal
deformation of the formal group law is classified by the Lubin--Tate ring
\[R:= W(\bb{F}) \llbracket v_1, \ldots, v_{n-1} \rrbracket \] 
where $W(\bb{F})$ is the ring of Witt vectors. There is an associated even
periodic cohomology theory $\e_n$, called the $n$-th Morava $E$-theory or
Lubin-Tate theory, with
\[\pi_*\e_n \cong W(\bb{F}) \llbracket v_1, \ldots, v_{n-1} \rrbracket [u^{\pm}],\] 
where the $v_i$'s are of degree zero and $u$ is of degree $-2$. Note that the
spectrum $\e_n$ depends not only on $n$ but also on the choice of $\bb{F}$ and
$\Gamma$. In this paper we will take $\mathbb F=\mathbb F_{p^{p-1}}$.

We are interested in higher real $K$-theories of height $p - 1$ at a prime $p$.
These theories are constructed as follows. Recall that the Morava $E$-theory $\e_{n}$ comes equipped with an action of the Morava
stabilizer group $\text{Aut}(\Gamma)$. 
When $n=p-1$  
there is a maximal subgroup $G$ of $\text{Aut}(\Gamma)$ containing
an element of order $p$, unique up to conjugacy. We define 
\[\eop= \e_{p\,\text{-}1}^{\text{h}G},\]
the homotopy fixed points of $\e_{p\,\text{-}1}$ with respect to $G$.
These spectra are called
higher real $K$-theories because $\eo_1$ at the prime $2$
is the $K(1)$-localization of the real $K$-theory. When $p=3$, $\eo_2$ is the
$K(2)$-localization of topological modular forms. 

The homotopy
groups of the higher real $K$-theories $\eop$ are computed by Hopkins and
Miller using the homotopy fixed point spectral sequence. See \cite[Section
1]{Chat} for an illustration with information about elements in the Hurewicz
image.  When the prime is understood, the spectrum $\eop$ will be abbreviated by $\eo$
 and $\e_{p\,\text{-}1}$ by $\e$.

\subsection{\texorpdfstring{$\eop$}{EO}-modules and algebraic \texorpdfstring{$\eop$}{EO}-theory}\label{subsec:alg_eo}

We now summarize results from \cite{Chat} describing the $\eo$-module theory of certain spectra. 
These are key technical tools we will use extensively.

\begin{defn}\label{defn:alg_eo_stuff} Let $p$ be a prime. 
\begin{enumerate}
\item Let $P(1)^*$ be the sub Hopf algebra of the mod 
$p$-Steenrod algebra generated by the power operation 
$P^1$. Let $P(1)_*\simeq \F_p[\xi_1]/\xi_1^p$ denote the 
dual quotient Hopf algebra of the dual Steenrod algebra. 
\item  Let $W_l$ denote the unique indecomposable length 
$l$ graded $P(1)_*$-comodule, concentrated in degrees $0$ through $2(l-1)(p-1)$.
\item Let $\Sigma^sW_l$ denote the indecomposable length 
$l$ graded $P(1)_*$-comodule concentrated in degrees $s$ through $s+2(l-1)(p-1)$.
\item 
Let $X_l$ denote the spectrum uniquely determined by the property that $H_*(X_l;\Fp) \simeq W_l$ 
as a $P(1)_*$-comodule and such that all Steenrod operations $P^i$ for $i>1$ act by zero.
\end{enumerate}
\end{defn}

\begin{defn}[{\cite[Definition 5.2]{Chat}}]
\label{def:algeo} 
    A spectrum $Z$ has {\em algebraic $\eo$-theory} if there exists an index set $I$, $s_i\in \mathbb Z$ and $l_i \in \{1,\ldots, p\}$ such that
\[
    \eo \otimes Z \simeq \eo \otimes \left( \oplus_{i\in I} \Sigma^{s_i}X_{l_i}\right)
\]
 in the category of $\eo$-modules.
\end{defn}

\begin{ex}\label{ex:cpnr} The following spectra have algebraic $\eo$-theory.
\begin{itemize}
\item The spectrum $\cp{r}{n}$. To see this, recall
that the $p$-completed $\cp{r}{n}$ has an Adams splitting \cite{McG}. Each Adams
summand has cells in degrees differing by $2p-2$. Thus $\cp{r}{n}$ is a sum of
$2p-2$-sparse spectra. These summands have algebraic $\eo$-theory by
\cite[Theorem 5.13]{Chat} and therefore so does $\cp{r}{n}.$
\item The spectrum $\cp{r}{n} \otimes D \cp{r}{n}$. This
follows from the previous item and the fact that the tensor product of algebraic
$\eo$-modules is algebraic \cite[Corollary 5.24]{Chat}.
\end{itemize}
\end{ex}
The $\eo$-module structure of a spectrum
 $Z$ with algebraic $\eo$-theory is
determined by the $P^1$-module structure on $\Fp$-cohomology. 

\begin{thm}[{\cite[Theorem 5.6]{Chat} }]     
Let $Z$ be a spectrum with algebraic $\eo$-theory. Suppose that
$H_*(Z)\simeq \oplus_i \Sigma^{s_i}W_{l_i}$ as a $P(1)_*$-comodule. Then $\eo
\otimes Z \simeq \eo\otimes\left( \oplus_i \Sigma^{s_i} X_{l_i}\right)$
in the category of $\eo$-modules.
\end{thm}

Another key technical tool from \cite{Chat} is the following orientability result.

\begin{thm}[{\cite[Corollary 1.6]{Chat}}]
    \label{thm:eo_orient} 
    Let $\gamma_1$ denote the universal bundle on $\cp{}{\infty}$. Then
    $p\gamma_1$ is $\eo$-orientable.
\end{thm}
This implies that $p$ times the tautological bundle $\gamma_1$ on $\cp{}{n}$ is
$\eo$-orientable. Given an integer $k$, let $[k]_p$ denote the remainder of $k$
modulo $p$ in $\{0,\ldots, p-1\}$.  
Since $\cp{r}{n}$ is isomorphic to the Thom spectrum $({\cp{}{n-r}})^{{r
\gamma_1}}$, \Cref{thm:eo_orient} implies that 
\[\eo \otimes \cp{r}{n} \simeq \Sigma^{2(r-[r]_p)}\eo \otimes \cp{[r]_p}{n-r+[r]_p}.\]

Together, \Cref{thm:eo_orient} and the previous observation imply:
\begin{cor}
\label{cor:eo_periodicity}
Given a fixed integer $c$, the group $\hmaps{\eo}{\eo \otimes \cp{r}{r+c}}{\eo \otimes \Sigma \cp{r}{r+c}}$ is
periodic in $r$ with period $p$.
\end{cor}

\Cref{cor:eo_periodicity} can be viewed as an improvement of James or Atiyah--Todd periodicity (cf. \Cref{rm:James}).  We make use of this smaller periodicity as
follows: the lower bounds for vector bundle counts are achieved by computing $p$ cases for any fixed
corank $c$ and inferring all by periodicity.

\section{First calculations}\label{sec:first}

Fix a prime $p$. Throughout the rest of the paper, all spectra are implicitly completed at $p$, and all homotopy groups and groups of stable vector bundles are understood to be up to $p$-completion.
In this section, we present direct calculations for $\hmapspp{\cp{r}{n}}{\Sigma\cp{r}{n}}$ and discuss its 
periodic behavior with respect to the rank $r$, when the corank $c:=n-r$ is relatively small (no greater than $2p-2$).
 Note that the metastable range translates to $0<c\leq r$.

We work with the cohomology of stunted projective spaces $\cp{a}{b}$, necessitating the next definition.
\begin{defn}\label{defn:xi} Let $x$ be the first Chern class of the tautological line bundle on $\cp{}{b}$. 
For $a \leq i \leq b$, write $x_i$ for the pullback of $x^i \in H^{2i}(\cp{}{b})$ along the quotient map $q: \cp{}{b} \rightarrow \cp{a}{b}$. 
\end{defn}

\subsection{Calculations at corank \texorpdfstring{$p-1$}{p−1}}\label{subsec:small_corank1}
To start, we note that
\[ \hmapspp{\cp{r}{r+c}}{\Sigma\cp{r}{r+c}} = 0\]
for $c<p-1$.
 Indeed, at an odd prime $p$ the first nontrivial stable homotopy group is $\pi_{2p-3}(\sphere)\cong \Z/p$ generated by $\alpha_1$, 
as in \cite[1.3.10 and 1.3.11]{Ravenel_green}. When $c<p-1$, 
the dimension difference between the top and bottom cells of 
$\cp{r}{r+c}$ is no greater than $2p-4$, and this forces $\cp{r}{r+c}$ 
to split as a sum of spheres. 
Now ${\hmapspp{\cp{r}{r+c} }{ \Sigma\cp{r}{r+c}} }$ is a sum of groups of the form 
$\hmapspp{\sphere^{2k}}{\sphere^{2l+1}}$ where $2k-(2l+1) \leq 2p-5$.
Since $2p-5$ is strictly less than the degree of $\alpha_1$, the first $p$-torsion in $\hmapspp{\cp{r}{r+c} }{ \Sigma\cp{r}{r+c}}$ 
cannot appear when the corank $c$ is less than $p-1$.

The first interesting $p$-complete calculation is that of \[\hmapspp{\cp{p}{2p-1}}{ \Sigma\cp{p}{2p-1}}.\]
 The only possibly nontrivial Steenrod operation in the mod $p$ cohomology of $\cp{p}{2p-1}$ is
$P^1$ on $x_p$, where we refer the reader to \Cref{defn:xi} 
for the definition of $x_p$ and the map $q$ that appears below. It follows from the Cartan formula that
\[P^1(x_p) = P^1(q^*x^p) = q^*P^1(x^p) = q^*(px^{2p-1}) = 0,\]
and hence the $p$-complete splitting
\[\cp{p}{2p-1} \simeq \sphere ^{2p} \oplus \sphere ^{2p+2} \oplus \cdots \oplus \sphere^{4p-2}.\] 
Therefore
\[
\hmapspp{\cp{p}{2p-1}}{ \Sigma\cp{p}{2p-1}} 
\cong \hmapspp{ \bigoplus_{i=0}^{p-1}\sphere^{2p+2i}}{\,\, \bigoplus_{j=0}^{p-1}\sphere^{2p+1+2j}}  \cong  \pi_{4p-2}(\sphere^{2p+1})
\]
where the last group is $ \bb{Z}/p$ generated by $\alpha_1$. For any $r$
divisible by $p$, the same argument proves that 
\[\hmapspp{\cp{r}{r+p-1}}{ \Sigma\cp{r}{r+p-1}} \cong \bb{Z}/p.\]

Now consider $r\equiv a$ mod $p$, where $0<a\leq p-1$. By the Cartan formula \[P^1(x_r) = ax_{r+p-1},\] 
which generates the top-degree cohomology of $\cp{r}{r+p-1}$. This implies a $p$-complete splitting
\[\cp{r}{r+p-1} \simeq \Sigma^{2r}C(\alpha_1) \oplus \sphere^{2r+2} \oplus \cdots \oplus \sphere^{2r + 2p-4},\]
where $C(\alpha_1)$ denotes the cone of $\alpha_1$.
So $\hmapspp{\cp{r}{r+p-1}}{ \Sigma\cp{r}{r+p-1}}$ breaks into a sum of abelian groups $\hmapspp{X}{ Y}$, where:
\begin{itemize}\item $X = C(\alpha_1)$ or some $\sphere^{2m}$ for $1\leq m \leq p-2$; and 
\item $Y = \Sigma C(\alpha_1)$ or some $\sphere^{2m+1}$ for $1\leq m \leq p-2$.
\end{itemize} When both $X$ and $Y$ are spheres, $\hmapsp{X}{ Y} = 0$ for dimension reasons.
 The group is also zero when $X$ is a sphere and $Y$ is $\Sigma C(\alpha_1)$, because $\pi_i(C(\alpha_1)) = 0$ for $0< i \leq 2p-5$. 
When $X = C(\alpha_1)$ and $Y$ is an odd sphere, we again have $\hmapsp{X}{Y} = 0$. Indeed, for any $1\leq m \leq p-2$, the group $\hmapspp{C(\alpha_1)}{\sphere^{2m+1}}$ is an extension of $\pi_0(\sphere^{2m+1})=0$ by $\pi_{2p-2}(\sphere^{2m+1})=0$. 
Finally, when $X = C(\alpha_1)$ and 
$Y = \Sigma C(\alpha_1)$, we get zero since the $\alpha_1$ from the top cell of $X$ to the bottom cell of $Y$ is killed by the $\alpha_1$-attaching map in $Y$. We have proved the following:

\begin{prop} The abelian group $\hmapspp{\cp{r}{r+p-1}}{\Sigma\cp{r}{r+p-1}}$ is isomorphic to $\bb{Z}/p$ if $r\equiv 0\pmod p$ and is zero
otherwise.
\end{prop}

This computation is a slightly more explicit warm-up for the next section.
\subsection{Calculations at corank \texorpdfstring{$c\leq 2p-3$}{2p−3}}\label{subsec:small_corank2}

To compute  $\hmapspp{\cp{r}{r+c}}{ \Sigma\cp{r}{r+c}}$, note that the spectrum $\cp{r}{r+c}$
 again splits as a sum of even suspensions of cofibers of $\alpha_1$ and even-dimensional spheres, with attaching maps determined by $r\pmod{p}$.

\begin{thm}\label{thm:count1}
Let $c$ be a fixed integer with $p-1 \leq c \leq 2p-3$. The abelian groups $\hmapsp{\cp{r}{r+c}}{ \Sigma\cp{r}{r+c}}$
 exhibit a $p$-fold periodicity in $r$. More precisely: 

\[ \hmapspp{\cp{r}{r+c}}{ \Sigma\cp{r}{r+c}}
 \cong \begin{cases} \bb{Z}/p & \text{ if some number in the list } \{r, \ldots, r+c-(p-1)\} \text{ is divisible by $p$}  \\
0 & \text{otherwise}. \end{cases}\] 
\end{thm}
\begin{proof}
    We break the computation into two cases. First suppose that there is some integer $k$ so that $r\leq k \leq r+c-(p-1)$, and that $k \equiv 0 \pmod p$. 
Note that \[P^1(x_i) = ix_{i+p-1},\]
and therefore $P^1(x_i)$ is nonzero for $i \in\{ r, r+1, \cdots, r+c-(p-1)\} \setminus \{k \}.$ On the other hand, $P^1(x_k)=0$. So there is a $p$-complete splitting (which in fact arises from the Adams splitting for $\cp{}{\infty}$):
    \begin{equation}\label{eq:decomposition-key}
\cp{r}{r+c} \simeq \left(
\bigoplus_{\substack{i=0 \\ r+i \neq k}}^{c-(p-1)} \Sigma^{2(r+i)}C(\alpha_1) \right)
 \oplus \sphere^{2k} \oplus \sphere^{2k+2p-2} \oplus X,\end{equation}
    where $X$ is a sum of spheres in each even dimension $j$ such that $2(r+c-(p-1))<j<2(r+p-1)$. It follows that
    \begin{equation}\label{identification-0}
\hmapspp{\cp{r}{r+c}}{ \Sigma \cp{r}{r+c}} \cong \pi_{2k+2p-3}(\bb{S}^{2k}) \cong \bb{Z}/p\cdot \alpha_1. \end{equation}
    
Second, suppose that there is no integer $k$ in the range $r\leq k \leq r+c-(p-1)$ so that $k \equiv 0 \pmod p$.
 In this case $P^1(x_i) = ix_{i+p-1}$ is nonzero for all $i = r, r+1, \cdots, r+c-(p-1)$. It then follows that there is a $p$-complete splitting
    \begin{equation}\label{adams-decomp-1}
\cp{r}{r+c} \simeq \left( \bigoplus_{i =0}^{c-(p-1)} \Sigma^{2(r+i)}C(\alpha_1) \right) \oplus X,\end{equation}
    where $X$ is a sum of spheres in each even dimension $j$ such that $2(r+c-(p-1))<j<2(r+p-1)$. It follows that 
    \[\hmapspp{\cp{r}{r+c}}{ \Sigma \cp{r}{r+c}} =0. \]

Finally, note that for a given $r$ and $c$ with $c\leq 2p-3$, there is an integer $k$ so that $r\leq k \leq r+c-(p-1)$
 and that $k \equiv 0 \pmod p$ if and only if $r \equiv 0, -1, \cdots, p-1-c  \pmod p$. This completes the proof.
\end{proof}

\begin{cor}\label{cor:bundle_ct1}
    Let $c$ be an integer so that $p-1\leq c \leq 2p-3$. In the metastable range, the $p$-component of the number of stably trivial rank $r$ bundles over $\cp{}{r+c}$ exhibits a $p$-fold periodicity. More precisely, it equals
    \begin{itemize}
        \item $p$, if $r+i\equiv 0 \pmod p$ for some $i\in \{ 0,\ldots,c-(p-1)\}$, and
        \item $1$, otherwise.
    \end{itemize}
\end{cor}

\subsection{A geometric interpretation of $p$-torsion for $p-1\leq c \leq 2p-3$.}\label{geometric-easy}
The proof of \Cref{thm:count1} actually gives some insight into the geometry of the bundles that we have discovered. After decomposing $ \cp{r}{r+c}$ into its Adams summands as in \Cref{eq:decomposition-key}, the non-trivial vector bundles detected arise from maps between a free cell 
$\sphere^{2k+2p-2}$ in the source $\cp{r}{r+c}$ and a cell $\Sigma \sphere^{2k}$ in the target 
$\Sigma \cp{r}{r+c}$.

If $r+c=k+p-1$ for $k\equiv 0 \pmod{p}$, the vector bundle $V$ detected can be constructed explicitly as follows: let 
$\sigma\: S^{2k+2p-2} \to BU(r)$ generate $p$-torsion in the stable homotopy groups of unitary groups, detected by $\alpha_1$ on the split $(2k+1)$-cell in the spectrum $\Sigma \cp{r}{r+c}$ under the identification $\pi_{2k+2p-2}BU(r) \cong \pi_{2k+2p-2}^{\op{st}}\Sigma \cp{r}{k+p-1}.$
A non-trivial, stably trivial bundle $V$ detected by \Cref{thm:count1} above is obtained by the composite
\begin{equation}\label[diagram]{diagram-V} \begin{tikzcd} 
\cp{}{k+p-1} \ar[r] \ar[rr, bend right, "V"]&S^{2k+2p-2} \ar[r,"\sigma"] &BU(r),
\end{tikzcd}
\end{equation}
where the first arrow is the quotient of $\cp{}{k+p-1}$ by the $(2k+2p-4)$-skeleton.  
Tracing through \Cref{eq:decomposition-key} and \Cref{identification-0}, this construction gives an explicit representative for the $p$-torsion bundles in $\hmapspp{\cp{r}{r+c}}{\Sigma \cp{r}{r+c}}$  for those $r$ and $c \leq 2p-3$ such that $r+c=k+p-1$, with $k\equiv 0 \pmod{p}$. 

If $k\equiv 0 \pmod{p}$ and $k+2p-2>r+c>k+p-1$, then the cell $\sphere^{2k+2p-2}$ in the Adams decomposition \Cref{eq:decomposition-key} is still the top cell in an Adams summand, although it is evidently not the top cell in $\cp{}{r+c}.$ Note that the passage from stably trivial vector bundles to stable maps to $\Sigma \cp{r}{\infty}$ is functorial in cohomologically even cell complexes. Let $W$ be a non-trivial, stably trivial vector bundle on $\cp{}{r+c}$ detected by the map \[\cp{r}{r+c}\to \sphere^{2k+2p-2} \xrightarrow{\alpha_1}\Sigma \sphere^{2k} \to \Sigma \cp{r}{r+c},\] where the first and last maps are projections onto or inclusions of the indicated summands from \Cref{eq:decomposition-key}.  Restricting along $\cp{}{k+p-1} \to \cp{}{r+c}$ and applying functoriality and the previous argument shows that we have a commutative diagram:

\[\begin{tikzcd} \cp{}{r+c} \ar[r,dashed, "W"] & BU(r) \\
\cp{}{k+p-1} \ar[u]\ar[ur, "\,\,\,V" below] &,
\end{tikzcd}
\]
where $V$ is as in \Cref{diagram-V}. Moreover, after $p$-completion, $W$ is the unique dashed arrow making the diagram commute. 

So, the vector bundle $V\: \cp{}{k+p-1}\to BU(r)$ as in \Cref{diagram-V} extends over $\cp{}{r+c}$ as long as $r+c<k+2p-2$, and the extension is unique, provided that we ask for an extension that is $p$-torsion in $\Vect_r^0(\cp{}{r+c})$. 
\begin{rmk}
Note that, if $r+c \geq k+2p-2$, then $V$ cannot extend, since there is a nontrivial cellular attaching map that obstructs an extension. In other words, $V$ extends (essentially uniquely) over higher-dimensional complex projective spaces exactly as far as one might hope it does. 
\end{rmk}

In light of this discussion, the result of \Cref{thm:count1} can be restated as follows: for coranks up to $2p-3$, the only $p$-torsion stably trivial vector bundles on complex projective spaces are: (1) those pulled back from spheres $S^{2k+2p-2}$ with $k\equiv 0 \pmod{p}$ under projections onto top cells and (2) the extensions of bundles of type (1) over cellular inclusions, which are $p$-locally unique.  This is \Cref{thm:Y}.


\section{Computations at corank \texorpdfstring{$c = 3,4$}{c=3,4}}
\label{subsec:corank2p2} 
In this section, we give some explicit examples to demonstrate how the enumeration of stably trivial vector bundles changes for $p$ small relative to the corank. In general, as $c$ grows, more elements in $\pi_*\sphere$ can contribute to \begin{equation}\label{keygroup}\hmapspp{\cp{r}{r+c}}{\Sigma \cp{r}{r+c}}.\end{equation} 
For example, at an odd prime $p$, the next nontrivial element
in the stable stem after $\alpha_1$ is $\alpha_2$, which generates a $\Z/p$ in
degree $4p-5$ \cite[1.3.11]{Ravenel_green}. The first opportunity for $\alpha_2$ to contribute
to the group \eqref{keygroup} occurs when $c=2p-2$. 
In this section, we complete the computation of corank $c=3$ stably trivial vector bundles on complex projective spaces, and also compute the corank $c=4$ case at the prime $p=3$.  

\subsection{An enumeration of corank 3 stably trivial vector bundles}\label{subsec:corank3}

Stably trivial corank $3$ bundles on $\cp{}{n}$ are metastable for $n$ at least $6$. So, for $r\geq 3$, $\Vect_r^0(\cp{}{r+3})$ is an abelian group. The group is $24$-periodic in $r$, and is a direct sum of $2$- and $3$-power torsion. For $p=3$, $2p-3\leq 3$, so the $3$-power torsion can be deduced from \Cref{thm:B}:
\begin{prop}\label{thm:corank3-3}
   Let $r\geq 3$. At $p=3$, the abelian group $\hmapspp{\cp{r}{r+3}}{\Sigma \cp{r}{r+3}}$ exhibits the following $3$-fold periodic behavior:
   \vspace{.5em}
    \begin{center}
        \begin{tabular}{|c|c|c|c|c|c|c|c|c|c|}
        \hline
            $r \mod 3$ & $0$ & $1$ & $2$  \\
        \hline
            $\hmapspp{\cp{r}{r+3}}{\Sigma \cp{r}{r+3}}$ 
& $\bb{Z}/3$ 
& $0$
& $\bb{Z}/3$ 
\\
\hline
\end{tabular}
\end{center}
\end{prop} 

The $2$-power torsion requires more information, since the prime $p=2$ is small relative to the corank $3$. We compute the orders of the groups in question, without solving more subtle extension problems.

\begin{prop} \label{thm:corank3-2}
   Let $r\geq 3$. At $p=2$, the abelian group $\hmapspp{\cp{r}{r+3}}{\Sigma \cp{r}{r+3}}$ exhibits the following $8$-fold periodic behavior:
   \vspace{.5em}
    \begin{center}
        \begin{tabular}{|c|c|c|c|c|c|c|c|c|c|}
        \hline
            $r \mod 8$ & $0$ & $1$ & $2$ & $3$ & $4$ & $5$ & $6$ & $7$ \\
        \hline
            $\hmapspp{\cp{r}{r+3}}{\Sigma \cp{r}{r+3}}$ 
& $\bb{Z}/4 \oplus \Z/2$ 
& $\Z/2?\Z/2 $ 
& $0$
& $\Z/4$ 
& $\Z/2? \Z/4$ 
& $\bb{Z}/4$ 
& $0$ 
& $\bb{Z}/8$ \\
        \hline
        \end{tabular}     
    \end{center}
  \vspace{.5em}
where $\Z/2^k \,?\, \Z/2^m$ indicates an unresolved extension problem, whose resolution is independent of $r\pmod{8}$.
\end{prop}

\begin{proof}
If $r$ is even, we can produce an explicit extension of $\hmapsp{\cp{r}{r+3}}{\Sigma \cp{r}{r+3}}$ by two known groups without extensive computation. If $r$ is odd, we need to do a bit more work.

First, it's useful to note the stable homotopy groups $\pi_{2r+i}\Sigma \cp{r}{r+3}$ for $1 \leq i \leq 6$, which are essentially classical and can be computed by the following methods:
 \begin{itemize}
\item $i=1,2,3$: $\pi_{2r+i}\Sigma\cp{r}{r+3}$ is isomorphic to $\pi_{i-1}C(\eta)$ for $i$ odd, and $\pi_{i-1}(\sphere^0 \oplus \sphere^2)$ for $i$ even.
\item $i=4,6$: these groups are isomorphic to $\pi_{2r+(i-1)}U(r)$, which are computed in \cite[page 971]{Mimura_HBAT};
\item $i=5$: these groups can be computed easily using long exact sequences. From the fiber sequence $\cp{r}{r+2} \to \cp{r}{r+3} \to \sphere^{2r+6},$ we find that \[\pi_{2r+4}\cp{r}{r+2} \cong \pi_{2r+4}\cp{r}{r+3}\cong \pi_{2r+5}\Sigma\cp{r}{r+3}.\] For $r$ even, from the fiber sequence $\sphere^{2r} \to \cp{r}{r+2} \to \cp{r+1}{r+2} $ we find that 
\[\pi_{2r+4}\cp{r}{r+2}\hookrightarrow \pi_{2r+4}\cp{r+1}{r+2}\cong \Z,\] using the identification $\cp{r+1}{r+2}\cong \Sigma^{2r+2}C(\eta)$. The multiple of the generator for $\pi_{2r+4}\Sigma^{2r+2}C(\eta)$ that is in the image is determined by a further attaching map.
For $r$ odd, the same fiber sequence gives an exact sequence
\begin{equation}\label[empty]{LES-26}0\to \pi_{2r+4}\cp{r}{r+2}\to \pi_{2r+4}\cp{r+1}{r+2}\to  \pi_{2r+4}\sphere^{2r+1}.\end{equation}
Note that $\pi_{2r+4}\cp{r+1}{r+2}\cong \Z/2 \oplus \Z$, using the identification $\cp{r+1}{r+2}\cong
\sphere^{2r+2}\oplus \sphere^{2r+4}$. The map $\cp{r+1}{r+2} \simeq \sphere^{2r+2} \oplus \sphere^{2r+4} \to \sphere^{2r+1}$ restricts to $\eta$ on the first summand, so that the kernel of the last map in \Cref{LES-26} includes the summand $\Z/2$. The multiple of the generator for $\pi_{2r+4}S^{2r+4}$ of the remaining summand that is actually in the image of the map from $\pi_{2r+4}\cp{r}{r+2}$ is determined by a further attaching map.
\end{itemize} 
The relevant $2$-complete homotopy of $\Sigma\cp{r}{r+3}$ is indicated in \Cref{stunted-htpy}.
\begin{figure}[h]
  \begin{tabular}{|c|c|c|c|c|c|c|c|}
    \hline
    $q$
      & $2r+1$
      & $2r+2$
      & $2r+3$
      & $2r+4$
      & $2r+5$
      & $2r+6$ \\
    \hline
    $\pi_q\Sigma\cp{r}{r+3},$ $r\equiv 0 \pmod{8}$
      & $\bb{Z}$
      & $\bb{Z}/2$
      & $\bb{Z}\oplus\bb{Z}/2$
      & $\bb{Z}/8$
      & $\bb{Z}$
      & $\bb{Z}/2$\\
    \hline
    $\pi_q\Sigma\cp{r}{r+3},$ $r \equiv 1 \pmod{8}$
      & $\bb{Z}_{}$
      & $0$
      & $\bb{Z}_{}$
      & $\bb{Z}/2$
      & $\bb{Z}_{}$
      & $\bb{Z}/2$ \\
    \hline
    $\pi_q\Sigma\cp{r}{r+3}$, $r \equiv 2 \pmod{8}$
      & $\bb{Z}$
      & $\bb{Z}/2$
      & $\bb{Z}\oplus\bb{Z}/2$
      & $\bb{Z}/2$
      & $\bb{Z}$
      & $0$ \\
    \hline
    $\pi_q\Sigma\cp{r}{r+3},$ $r\equiv 3 \pmod{8}$
      & $\bb{Z}$
      & $0$
      & $\bb{Z}$
      & $0$
      & $\bb{Z}$
      & $\bb{Z}/4$
\\
    \hline
 $\pi_q\Sigma\cp{r}{r+3},$ $r \equiv 4 \pmod{8}$
      & $\bb{Z}$
      & $\bb{Z}/2$
      & $\bb{Z}\oplus\bb{Z}/2$
      & $\bb{Z}/4$
      & $\bb{Z}$
      & $\bb{Z}/4$
\\
    \hline
    $\pi_q\Sigma\cp{r}{r+3},$ $r \equiv 5 \pmod{8}$
      & $\bb{Z}$
      & $0$
      & $\bb{Z}$
      & $\bb{Z}/4 $
      & $\bb{Z}$
      & $\bb{Z}/2$\\
    \hline
 $\pi_q\Sigma\cp{r}{r+3}$, $r \equiv 6 \pmod{8}$
      & $\bb{Z}$
      & $\bb{Z}/2$
      & $\bb{Z}\oplus\bb{Z}/2$
      & $\bb{Z}/2$
      & $\bb{Z}$
      & $0$ \\
    \hline
    $\pi_q\Sigma\cp{r}{r+3},$ $r\equiv 7\pmod{8}$
      & $\bb{Z}$
      & $0$
      & $\bb{Z}$
      & $0$
      & $\bb{Z}$
      & $\bb{Z}/8$\\
    \hline
  \end{tabular}
\caption{The low-degree stable homotopy of $\Sigma \cp{r}{r+3}$, up to $2$-completion.}\label{stunted-htpy}
\end{figure}

Now consider the fiber sequence of spectra
\begin{equation}\label{useful-1} \cp{r}{r+2}\to \cp{r}{r+3} \to \sphere^{2r+6}.\end{equation}
Applying $\hmapsp{-}{\Sigma \cp{r}{r+3}}$, we claim that we get an exact sequence:
\begin{equation}\label[empty]{LES-11}
\hmapsp{\sphere^{2r+6}}{\Sigma \cp{r}{r+3}} \to \hmapsp{\cp{r}{r+3}}{\Sigma \cp{r}{r+3}} \to \hmapsp{\cp{r}{r+2}}{\Sigma \cp{r}{r+3}} \to 0.
\end{equation}
The content is to prove that the map \[\hmapsp{\cp{r}{r+3}}{\Sigma \cp{r}{r+3}} \to \hmapsp{\cp{r}{r+2}}{\Sigma \cp{r}{r+3}}\] is surjective. Note that the next term to the right in the long exact sequence is \[\hmapsp{\sphere^{2r+5}}{\Sigma \cp{r}{r+3}}\cong\pi_{2r+5}\Sigma\cp{r}{r+3} \cong \Z\] by \Cref{stunted-htpy}. On the other hand, 
 \[\hmapsp{\cp{r}{r+2}}{\Sigma \cp{r}{r+3}} \cong \hmapsp{\cp{r}{r+2}}{\Sigma \cp{r}{r+2}}\cong \Vect_r^0(\cp{}{r+2})\] is a torsion group, identified in \cite[Theorem 1.2]{Hu}.

We next study the first map
$
\hmapsp{\sphere^{2r+6}}{\Sigma \cp{r}{r+3}} \to \hmapsp{\cp{r}{r+3}}{\Sigma \cp{r}{r+3}} $ in \Cref{LES-11}.
Under the identifications \[\Vect_r^0(S^{2r+6}) \cong \hmapsp{\sphere^{2r+6}}{\Sigma \cp{r}{r+3}}, \quad \Vect_r^0(\cp{}{r+3})\cong \hmapsp{\cp{r}{r+3}}{\Sigma \cp{r}{r+3}},\]
this map is induced by pullback along projection onto the top cell $\cp{}{r+3} \to S^{2r+6}$. It is useful to observe the following commutative diagram:

\begin{equation}\label{eq:comm-vary}
\begin{tikzcd} 
\pi_{2r+6}BU(r) \ar[r,"p^*"]\ar[d,"\cong"] & \Vect_r^0(\cp{}{r+3})\ar[d]\\
\pi_{2r+6}BU(r+1) \ar[r, "q^*"] & \Vect_{r+1}^0(\cp{}{r+3}),
\end{tikzcd}
\end{equation}
where, by \cite[page 971]{Mimura_HBAT}, the left vertical map is an isomorphism. By \cite[Theorem 1.3]{opie24enum}, any two rank $r+1$ bundles on $\cp{}{r+3}$ with the same $K$-theory class have the same restriction to $\cp{}{r+2}$. Hence, any stably trivial vector bundle restricts to zero on $\cp{}{r+2}$, and is pulled back from a rank $r+1$ bundle on $S^{2r+6}$ under the natural projection map to the top cell. In particular, the map $q^*$ is surjective. 

Comparing \cite[Theorem 1.2]{Hu} and the homotopy groups $\pi_{2r+6}BU(r+1)$ from \cite[Page 971]{Mimura_HBAT} or \Cref{stunted-htpy} above, we find that $q^*$ is bijective when $r+1$ is odd, and has kernel $\Z/2$ when $r+1$ is even.  

We now consider cases.
\begin{itemize}
\item[]
\begin{itemize}
\item[{\bf $r$ even}.] Since $q^*$ is bijective, the top horizontal map in \Cref{eq:comm-vary} has trivial kernel. Referring back to \Cref{LES-11}, we get a short exact sequence
\begin{equation}\label[empty]{even-extension} 0 \to \pi_{2r+6}BU(r) \to \Vect_r^0(\cp{}{r+3}) \to \Vect_r^0(\cp{}{r+2}) \to 0.\end{equation}
We can now compute $\Vect_r^0(\cp{}{r+3})$ up to extension, since $\pi_{2r+6}BU(r) \cong \pi_{2r+6}\Sigma \cp{r}{r+3}$ is given in \Cref{stunted-htpy}, while the group $\Vect^0_r(\cp{}{r+2})$ is in \cite[Theorem 4.2]{Hu}. To be explicit:
\[ \Vect_r^0(\cp{}{r+2})=\begin{cases} \Z/4, & r \equiv 0 \pmod{8}
\\
0, & r \equiv 2 \pmod{8} \\
\Z/2, & r \equiv 4 \pmod{8}\\
0, & r \equiv 6 \pmod{8}.
\end{cases}
\] 
In the case that $r\equiv 0 \pmod{8}$, we can actually see that the extension \Cref{even-extension} splits by noting that $\cp{r}{r+3} \cong \sphere^{2r} \oplus \cp{r+1}{r+3}$. 
\item[{\bf $r$ odd}.] From \Cref{LES-11} we get an exact sequence
\begin{equation}\label[empty]{eq:LES-odd} \Z/2 \to \pi_{2r+6}BU(r) \to \Vect_r^0(\cp{}{r+3}) \xrightarrow{a} \Vect_r^0(\cp{}{r+2}) \to 0.\end{equation}

More work is required. 
\end{itemize}
\end{itemize}

We now deal with each odd $r\pmod 8$ in cases. We compute portions of the Atiyah--Hirzebruch spectral sequence converging to stable maps from $\cp{r}{r+3}$ to $\Sigma \cp{r}{r+3}$. It is convenient to introduce the following notation for the Atiyah--Hirzebruch spectral sequence computing (shifted) stable maps. For any spectra $X,C$, where $X$ is finite, we have a strongly convergent spectral sequence:
\begin{equation}\label[diagram]{AHSS}\op{AHSS}(X,C):
  E_2^{q,p}
  =H^p(X;\pi_{-q}C)
  \implies
\pi_{-(p+q)}\mapsp{X}{C}.
\end{equation}
We first note a few useful observations:
\begin{rmk}\label{rmk:cellular-functoriality}
On each page of $\op{AHSS}(\cp{r}{r+3},\Sigma \cp{r}{r+3})$, one can project onto the bi-graded summand of terms with $p\leq 2r+4$. This is precisely implementing the map of Atiyah--Hirzebruch spectral sequences \[\op{AHSS}(\cp{r}{r+3},\Sigma \cp{r}{r+3})\to \op{AHSS}(\cp{r}{r+2},\Sigma \cp{r}{r+3})\] induced by cellular inclusion $\cp{r}{r+2} \to \cp{r}{r+3}$. We know this map is surjective on the $p+q=0$ line on $E_\infty$ pages, from surjectivity of the map labelled $a$ in \Cref{LES-11}.
\end{rmk} 

\begin{rmk}\label{rmk:permanent-cycle} We note that, independent of $r \pmod{8}$, in any spectral sequence $\op{AHSS}(\cp{r}{r+3},\Sigma \cp{r}{r+3})$ there is a $\Z$ term in bidegree $(-(2r+1), 2r)$ on the $E_2$-page. Recall that line $p+q=-1$ computes $\hmapsp{\Sigma \cp{r}{r+3}}{\Sigma \cp{r}{r+3}}$; moreover, the term in bidegree $(-(2r+1),2r)$ corresponds to the bottom cell map $\sphere^{2r+1} \to \Sigma \cp{r}{r+3}$. Whether or not it survives the spectral sequence amounts to extending the map over the cellular inclusion $\sphere^{2r+1} \to \Sigma \cp{r}{r+3}$. However, it extends to a map $\Sigma \cp{r}{r+3} \to \Sigma \cp{r}{r+3}$ via the identity map, so this term is always a permanent cycle and supports no nontrivial differentials.
\end{rmk}

\begin{conv}
In what follows, it is convenient to label generators of $\pi_{2r+i}{\Sigma\cp{r}{r+3}}$, $1\leq i\leq 6$, by appropriate elements in the stable stem to show where they come from. We write $\pi_{2r+i}{\Sigma\cp{r}{r+3}} \cong A\cdot \alpha[N]$, if the homotopy group is isomorphic to the abelian group $A$, whose generator is supported by $\alpha\in \pi_*\sphere^0$ on the standard class $\iota_N \in \pi_N{\sphere^N}$.
\end{conv}

\begin{itemize}
\item[{\bf Case \!\! 1}]: $r\equiv 1 \pmod{8}$. The cell structure of $\cp{r}{r+3}$ is depicted in \Cref{fig:r1mod8}. The attaching maps can be deduced from \cite[Table 2.4]{mosher69}.
\begin{equation}\label[diagram]{fig:r1mod8}
  \centering
  \begin{tikzpicture}[x=0.9cm,y=0.9cm,baseline=(current bounding box.center)]
    \fill (0,0) circle (2pt);
    \fill (2,0) circle (2pt);
    \fill (4,0) circle (2pt);
    \fill (6,0) circle (2pt);

    \draw[line width=0.5pt] (0,0) -- (2,0);
    \node[below] at (1,0) {$\eta$};

    \draw[line width=0.5pt] (4,0) -- (6,0);
    \node[below] at (5,0) {$\eta$};

\draw[line width=0.5pt]
  (6,0) to[out=155,in=25,looseness=1.15] (2,0);
\node[above=-5mm] at (3.15,1.15) {$\nu$};

    \draw[line width=0.5pt]
      (4,0) to[out=-35,in=-145,looseness=1.15] (0,0);
    \node[below=-2mm] at (2,-1.15) {$2\nu$};

    \node[below=16mm] at (0,0) {$2r$};
    \node[below=16mm] at (2,0) {$2r+2$};
    \node[below=16mm] at (4,0) {$2r+4$};
    \node[below=16mm] at (6,0) {$2r+6$};
  \end{tikzpicture}
\end{equation}
The relevant portion of the $E_2$-page for $\op{AHSS}(\cp{r}{r+3},\Sigma \cp{r}{r+3})$ is depicted in \Cref{E2-1}.  The generator in bidegree $(-(2r+4),2r+4)$ is a non-zero permanent cycle: consider \Cref{rmk:cellular-functoriality}. The $p+q=0$ line of $\op{AHSS}(\cp{r}{r+2},\Sigma \cp{r}{r+3})$ is identified with $\Vect_r^0(\cp{}{r+2})\cong \Z/2$, which forces the bidegree $(-(2r+4),2r+4)$ term to survive both spectral sequences.

It remains only to analyze the differentials with target the term $\Z/2$ in bidegree $(-(2r+6),2r+6)$. This group corresponds to restrictions of an element in $\pi_{2r+6}\Sigma \cp{r}{r+3}$ along projection onto the top cell $\cp{r}{r+3} \to \sphere^{2r+6}$. By \Cref{rmk:permanent-cycle}, there are two possible differentials with target this group: a $d_2$ and $d_4$. 
\begin{itemize}
\item[-] The $d_2$-differentials from the $(2r+4)$-column to the $(2r+6)$-column are induced by multiplication by $\eta$. The $d_2$-differential with target in bidegree $(-(2r+6), 2r+6)$ is zero: we note that  $\pi_{2r+5}(\Sigma\cp{r}{r+3}) \cong \Z\cdot 2[2r+5]$, and multiplying by $\eta$ gives zero in $\pi_{2r+6}(\Sigma\cp{r}{r+3})\cong \Z/2\cdot \eta[2r+5]$ after being multiplied by $\eta$. So both copies of $\Z/2$ survive to the $E_4$-page.

\item[-] The $d_4$-differentials from the $(2r+2)$-column to the $(2r+6)$ column are induced by multiplication by $\nu$ on homotopy.
The product of $\nu$ and the generator of $\pi_{2r+3}(\Sigma\cp{r}{r+3})\cong \Z\cdot 2[2r+3]$ does not survive in $\pi_{2r+6}(\Sigma\cp{r}{r+3}) \cong \Z/2 \cdot \eta[2r+5]$. So the $d_4$-differential from the group $\Z$ in bidegree $(-(2r+3), 2r+2)$ to the group $\Z/2$ in bidegree $(-(2r+6), 2r+6)$ is zero.
\end{itemize}
In conclusion, the $p+q=0$ line on the $E_\infty$-page is $\Z/2\oplus \Z/2$. We do not solve the extension problem.

\begin{figure}[h]
\centering
  \begin{tabular}{|c|c|c|c|c|c|c|c|}
    \hline
    $q\backslash p$
      & $2r$
      & $2r+1$
      & $2r+2$
      & $2r+3$
      & $2r+4$
      & $2r+5$
      & $2r+6$ \\
    \hline
    $-(2r+1)$
      & $\bb{Z}$ & $0$ & $\bb{Z}$ & $0$
      & $\bb{Z}$ & $0$ & $\bb{Z}$ \\
    \hline
    $-(2r+2)$
      & $0$ & $0$ & $0$ & $0$
      & $0$ & $0$ & $0$ \\
    \hline
    $-(2r+3)$
      & $\bb{Z}$ & $0$ & $\bb{Z}$ & $0$
      & $\bb{Z}$ & $0$ & $\bb{Z}$ \\
    \hline
    $-(2r+4)$
      & $\bb{Z}/2$ & $0$ & $\bb{Z}/2$ & $0$
      & $\bb{Z}/2$ & $0$ & $\bb{Z}/2$ \\
    \hline
    $-(2r+5)$
      & $\bb{Z}$ & $0$ & $\bb{Z}$ & $0$
      & $\bb{Z}$ & $0$ & $\bb{Z}$ \\
    \hline
    $-(2r+6)$
      & $\bb{Z}/2$ & $0$ & $\bb{Z}/2$ & $0$
      & $\bb{Z}/2$ & $0$ & $\bb{Z}/2$ \\
    \hline
  \end{tabular}
\caption{The $E_2$-page of the Atiyah--Hirzebruch spectral sequence computing $\pi_*\mapsp{\cp{r}{r+3}}{\Sigma \cp{r}{r+3}}$ when $r\equiv 1\pmod{8}$.}
\label{E2-1}
\end{figure}

\item[{\bf Case \!\! 2}]: $r\equiv 3 \pmod{8}$. The cell structure of $\cp{r}{r+3}$ is depicted in \Cref{fig:r3mod8}. The attaching maps follow from \cite[Table 2.4]{mosher69}.
\begin{equation}\label[diagram]{fig:r3mod8}
  \centering
  \begin{tikzpicture}[x=0.9cm,y=0.9cm,baseline=(current bounding box.center)]
    \fill (0,0) circle (2pt);
    \fill (2,0) circle (2pt);
    \fill (4,0) circle (2pt);
    \fill (6,0) circle (2pt);

    \draw[line width=0.5pt] (0,0) -- (2,0);
    \node[below] at (1,0) {$\eta$};

    \draw[line width=0.5pt] (4,0) -- (6,0);
    \node[above] at (5,0) {$\eta$};

    \draw[line width=0.5pt]
      (4,0) to[out=155,in=25,looseness=1.15] (0,0);
    \node[above=-4mm] at (2,1.15) {$\nu$};

    \draw[line width=0.5pt]
      (6,0) to[out=-155,in=-25,looseness=1.15] (2,0);
    \node[below=-4mm] at (4,-1.15) {$2\nu$};

    \node[below=16mm] at (0,0) {$2r$};
    \node[below=16mm] at (2,0) {$2r+2$};
    \node[below=16mm] at (4,0) {$2r+4$};
    \node[below=16mm] at (6,0) {$2r+6$};
  \end{tikzpicture}
\end{equation}
The relevant part of the $E_2$-page for $\op{AHSS}(\cp{r}{r+3},\Sigma \cp{r}{r+3})$ is in \Cref{E2-3}. 

We need to analyze the differentials with target the term $\Z/4$ in bidegree $(-(2r+6),2r+6)$. This group corresponds precisely to restrictions of an element in $\pi_{2r+6}\Sigma \cp{r}{r+3}$ along projection onto the top cell $\cp{r}{r+3} \to \sphere^{2r+6}$. By \Cref{rmk:permanent-cycle}, there are two possible differentials with target this group: a $d_2$ and $d_4$. 

\begin{itemize}
\item[-] The $d_2$-differential from the group $\Z$ in bidegree $(-(2r+5),2r+4)$ to the group $\Z/4$ in bidegree $(-(2r+6),2r+6)$ is induced by $\eta$. However, $\eta$ times the generator of $\pi_{2r+5}(\Sigma\cp{r}{r+3}) \cong \Z \cdot 8[2r+5]$ becomes zero in $\pi_{2r+6}(\Sigma\cp{r}{r+3}) \cong \Z/4 \cdot \nu[2r+3]$. So this $d_2$-differential is zero. 
\item[-] Consider the $d_4$-differential from the group $\Z$ in bidegree $(-(2r+3),2r+2)$ to the group $\Z/4$ in bidegree $(-(2r+6),2r+6)$. Note that it is a homomorphism induced by $2\nu$, from $\pi_{2r+3}(\Sigma\cp{r}{r+3}) \cong \Z \cdot 2[2r+3]$ to $\pi_{2r+6}(\Sigma\cp{r}{r+3}) \cong \Z/4 \cdot \nu[2r+3]$, so is also zero.
The $p+q=0$ line on the $E_\infty$-page therefore consists of a single summand $\Z/4$, and there are no extension problems.
\end{itemize}
\begin{figure}[h]
\centering
  \begin{tabular}{|c|c|c|c|c|c|c|c|}
    \hline
    $q\backslash p$
      & $2r$
      & $2r+1$
      & $2r+2$
      & $2r+3$
      & $2r+4$
      & $2r+5$
      & $2r+6$ \\
    \hline
    $-(2r+1)$
      & $\bb{Z}$ & $0$ & $\bb{Z}$ & $0$
      & $\bb{Z}$ & $0$ & $\bb{Z}$ \\
    \hline
    $-(2r+2)$
      & $0$ & $0$ & $0$ & $0$
      & $0$ & $0$ & $0$ \\
    \hline
    $-(2r+3)$
      & $\bb{Z}$ & $0$ & $\bb{Z}$ & $0$
      & $\bb{Z}$ & $0$ & $\bb{Z}$ \\
    \hline
    $-(2r+4)$
      & $0$ & $0$ & $0$ & $0$
      & $0$ & $0$ & $0$ \\
    \hline
    $-(2r+5)$
      & $\bb{Z}$ & $0$ & $\bb{Z}$ & $0$
      & $\bb{Z}$ & $0$ & $\bb{Z}$ \\
    \hline
    $-(2r+6)$
      & $\bb{Z}/4$ & $0$ & $\bb{Z}/4$ & $0$
      & $\bb{Z}/4$ & $0$ & $\bb{Z}/4$ \\
    \hline
  \end{tabular}
\caption{The $E_2$-page of the Atiyah--Hirzebruch spectral sequence when $r\equiv 3\pmod{8}$.}
\label{E2-3}
\end{figure}

\item[{\bf Case 3}]: $r\equiv 5 \pmod{8}$. The cell structure of $\cp{r}{r+3}$ is depicted in \Cref{fig:r5mod8}. The attaching maps can be deduced from \cite[Table 2.4]{mosher69}.
\begin{equation}\label[diagram]{fig:r5mod8}
  \centering
  \begin{tikzpicture}[x=0.9cm,y=0.9cm,baseline=(current bounding box.center)]
    \fill (0,0) circle (2pt);
    \fill (2,0) circle (2pt);
    \fill (4,0) circle (2pt);
    \fill (6,0) circle (2pt);

    \draw[line width=0.5pt] (0,0) -- (2,0);
    \node[below] at (1,0) {$\eta$};

    \draw[line width=0.5pt] (4,0) -- (6,0);
    \node[below] at (5,0) {$\eta$};

    \draw[line width=0.5pt]
      (6,0) to[out=155,in=25,looseness=1.15] (2,0);
    \node[above=-5mm] at (3.15,1.15) {$\nu$};

    \node[below=9mm] at (0,0) {$2r$};
    \node[below=9mm] at (2,0) {$2r+2$};
    \node[below=9mm] at (4,0) {$2r+4$};
    \node[below=9mm] at (6,0) {$2r+6$};
  \end{tikzpicture}
\end{equation}
The relevant portion of the $E_2$-page for $\op{AHSS}(\cp{r}{r+3},\Sigma \cp{r}{r+3})$ is depicted in \Cref{E2-5}. The generator in bidegree $(-(2r+4),2r+4)$ is a non-zero permanent cycle: by \Cref{rmk:cellular-functoriality}, the $p+q=0$ line of $\op{AHSS}(\cp{r}{r+2},\Sigma \cp{r}{r+3})$ is identified with $\Vect_r^0(\cp{}{r+2})\cong \Z/4$, which forces the bidegree $(-(2r+4),2r+4)$ term to survive both spectral sequences.

Analyzing cell structure and cohomology, there is a surjective $d_2$ differential with target the copy of $\Z/2$ in bidegree $(-(2r+6),2r+6)$. Therefore the group is $\Z/4$ and there are no extension problems.

\begin{figure}[h]
\centering
  \begin{tabular}{|c|c|c|c|c|c|c|c|}
    \hline
    $q\backslash p$
      & $2r$
      & $2r+1$
      & $2r+2$
      & $2r+3$
      & $2r+4$
      & $2r+5$
      & $2r+6$ \\
    \hline
    $-(2r+1)$
      & $\bb{Z}$ & $0$ & $\bb{Z}$ & $0$
      & $\bb{Z}$ & $0$ & $\bb{Z}$ \\
    \hline
    $-(2r+2)$
      & $0$ & $0$ & $0$ & $0$
      & $0$ & $0$ & $0$ \\
    \hline
    $-(2r+3)$
      & $\bb{Z}$ & $0$ & $\bb{Z}$ & $0$
      & $\bb{Z}$ & $0$ & $\bb{Z}$ \\
    \hline
    $-(2r+4)$
      & $\bb{Z}/4$ & $0$ & $\bb{Z}/4$ & $0$
      & $\bb{Z}/4$ & $0$ & $\bb{Z}/4$ \\
    \hline
    $-(2r+5)$
      & $\bb{Z}$ & $0$ & $\bb{Z}$ & $0$
      & $\bb{Z}$ & $0$ & $\bb{Z}$ \\
    \hline
    $-(2r+6)$
      & $\bb{Z}/2$ & $0$ & $\bb{Z}/2$ & $0$
      & $\bb{Z}/2$ & $0$ & $\bb{Z}/2$ \\
    \hline
  \end{tabular}
\caption{The $E_2$-page of the Atiyah--Hirzebruch spectral sequence when $r\equiv 5\pmod{8}$.}
\label{E2-5}
\end{figure}

\item[{\bf Case \!\! 4}]:   $r\equiv 7 \pmod{8}$. The cell structure of $\cp{r}{r+3}$ is depicted in \Cref{fig:r7mod8}. The attaching maps follow from \cite[Table 2.4]{mosher69}.

\begin{equation}\label[diagram]{fig:r7mod8}
  \centering
  \begin{tikzpicture}[x=0.9cm,y=0.9cm,baseline=(current bounding box.center)]
    \fill (0,0) circle (2pt);
    \fill (2,0) circle (2pt);
    \fill (4,0) circle (2pt);
    \fill (6,0) circle (2pt);

    \draw[line width=0.5pt] (0,0) -- (2,0);
    \node[below] at (1,0) {$\eta$};

    \draw[line width=0.5pt] (4,0) -- (6,0);
    \node[below] at (5,0) {$\eta$};

    \draw[line width=0.5pt]
      (0,0) to[out=25,in=155,looseness=1.15] (4,0);
    \node[above=-5mm] at (2,1.15) {$\nu$};

    \node[below=9mm] at (0,0) {$2r$};
    \node[below=9mm] at (2,0) {$2r+2$};
    \node[below=9mm] at (4,0) {$2r+4$};
    \node[below=9mm] at (6,0) {$2r+6$};
  \end{tikzpicture}
\end{equation}
The relevant part of the $E_2$-page for $\op{AHSS}(\cp{r}{r+3},\Sigma \cp{r}{r+3})$ is in \Cref{E2-7}. 

We need to analyze the differentials with target the term $\Z/8$ in bidegree $(-(2r+6),2r+6)$. This group corresponds precisely to restrictions of an element in $\pi_{2r+6}\Sigma \cp{r}{r+3}$ along projection onto the top cell $\cp{r}{r+3} \to \sphere^{2r+6}$. By \Cref{rmk:permanent-cycle}, there are two possible differentials with target this group: a $d_2$ and $d_4$. 

\begin{itemize}
\item[-] The $d_2$-differential from the $(2r+4)$-column to the $(2r+6)$-column is induced by $\eta$. Consider the one from the group $\Z$ in bidegree $(-(2r+5), 2r+4)$ to the group $\Z/8$ in bidegree $(-(2r+6), 2r+6)$. Since this is multiplication by $\eta$ from $\pi_{2r+5}\Sigma\cp{r}{r+3} \cong \Z\cdot 8[2r+5]$ to $\pi_{2r+6}\Sigma\cp{r}{r+3} \cong \Z/8\cdot \nu[2r+3]$, it must be the zero homomorphism.
\item[-] The $d_4$-differential from the $(2r+2)$-column to the $(2r+6)$-column is zero, since the cellular attaching map is zero.
\end{itemize} 

Therefore the only term on the $p+q=0$ line of the $E_\infty$-page is a copy of $\Z/8$.
\begin{figure}[h]
\centering
  \begin{tabular}{|c|c|c|c|c|c|c|c|}
    \hline
    $q\backslash p$
      & $2r$
      & $2r+1$
      & $2r+2$
      & $2r+3$
      & $2r+4$
      & $2r+5$
      & $2r+6$ \\
    \hline
    $-(2r+1)$
      & $\bb{Z}$ & $0$ & $\bb{Z}$ & $0$
      & $\bb{Z}$ & $0$ & $\bb{Z}$ \\
    \hline
    $-(2r+2)$
      & $0$ & $0$ & $0$ & $0$
      & $0$ & $0$ & $0$ \\
    \hline
    $-(2r+3)$
      & $\bb{Z}$ & $0$ & $\bb{Z}$ & $0$
      & $\bb{Z}$ & $0$ & $\bb{Z}$ \\
    \hline
    $-(2r+4)$
      & $0$ & $0$ & $0$ & $0$
      & $0$ & $0$ & $0$ \\
    \hline
    $-(2r+5)$
      & $\bb{Z}$ & $0$ & $\bb{Z}$ & $0$
      & $\bb{Z}$ & $0$ & $\bb{Z}$ \\
    \hline
    $-(2r+6)$
      & $\bb{Z}/8$ & $0$ & $\bb{Z}/8$ & $0$
      & $\bb{Z}/8$ & $0$ & $\bb{Z}/8$ \\
    \hline
  \end{tabular}
\caption{The $E_2$-page of the Atiyah--Hirzebruch spectral sequence when $r\equiv 7\pmod{8}$.}
\label{E2-7}
\end{figure}
\end{itemize}
\end{proof}

\subsection{Corank $4$ stably trivial vector bundles at the prime $3$}\label{subsec:corank4}

The enumeration of metastable, stably trivial corank $4$ vector bundles on complex projective spaces involves $p$-power torsion for $p=2,3,5$. The $2$-primary computation is beyond the scope of our examples here: it is $2^6$ periodic in $r$. The $5$-power torsion is covered by \Cref{thm:B}. The $3$-power torsion is $9$-periodic, and is not covered by \Cref{thm:B}. We work it out here. 

This provides some insight into how the odd-primary enumeration problem becomes more complicated when $c=2p-2$, which we can compare with partial information obtained in \Cref{cor:p3-case} below.

\begin{prop} \label{thm:corank_2p-2}
   At $p=3$, the abelian group $\hmapspp{\cp{r}{r+4}}{\Sigma \cp{r}{r+4}}$ exhibits a $9$-fold periodic behavior:
   \vspace{.5em}
    \begin{center}
        \begin{tabular}{|c|c|c|c|c|c|c|c|c|c|}
        \hline
            $r \mod 9$ & $0$ & $1$ & $2$ & $3$ & $4$ & $5$ & $6$ & $7$ & $8$ \\
        \hline
            $\hmapspp{\cp{r}{r+4}}{\Sigma \cp{r}{r+4}}$ & $\bb{Z}/9$ & $\bb{Z}/9$ & $\bb{Z}/3$ & $\bb{Z}/3$ 
& $\bb{Z}/3$ & $\bb{Z}/3$ & $\bb{Z}/3$ & $\bb{Z}/3$ & $\bb{Z}/3$\\
        \hline
        \end{tabular}.
    \end{center}
  \vspace{.5em}
\end{prop}

Since \Cref{thm:corank_2p-2} is meant to illustrate the challenges in extending direct computations beyond
 $c\leq 2p-3$, but is not necessary for what follows, we present the proof as an informal discussion.

The periodicity in the above calculation comes from that of the attaching maps in the cell structure of the five-cell stunted projective space $\cp{r}{r+4}$. When $r\equiv  0\pmod 9$, $\cp{r}{r+4}$ has the following cell structure
\begin{center}
\begin{tikzpicture}[scale = 0.8]
\fill (0,0) circle (2pt);
\fill (1,0) circle (2pt);
\fill (2,0) circle (2pt);
\fill (3,0) circle (2pt);
\fill (4,0) circle (2pt);
\draw (1,0) to [bend left] (3,0);
\draw (2,0) to [bend right] (4,0);
\end{tikzpicture}
\end{center}
where the bottom cell is of dimension $2r$ and the curved line segments indicate an $\alpha_1$-attaching map.
 In other words, $\cp{r}{r+4}$ splits as $\Sigma^{2r}\left(\sphere \oplus \Sigma^2C(\alpha_1) \oplus \Sigma^4C(\alpha_1)\right).$
 Therefore:
\[ \hmapspp{\cp{r}{r+4}}{ \Sigma \cp{r}{r+4}} \cong \hmapspp{\sphere\! \oplus\! \Sigma^2C(\alpha_1) \!\oplus\! \Sigma^4C(\alpha_1)}{ \sphere^1 \!\oplus\! \Sigma^3C(\alpha_1)\! \oplus\! \Sigma^5C(\alpha_1)}. \] 
Note that $\pi_0\left( \sphere^1 \oplus \Sigma^3C(\alpha_1) 
\oplus \Sigma^5C(\alpha_1)\right) =0$. Additionally, $\pi_*\left(C(\alpha_1)\right)$
 is zero in degree $2p-3 = 3$ and it is $\bb{Z}/9$ in degree 
$4p-5 = 7$. The latter comes from the extension
\[ 0 \rightarrow \pi_7(\sphere) \rightarrow \pi_7(C(\alpha_1)) \rightarrow \pi_7(\bb{S}^4) \rightarrow 0\]
whose nontriviality is witnessed by the Toda bracket $\langle p, \alpha_1, \alpha_1 \rangle = \alpha_2$. It follows that
\[ \hmapspp{\Sigma^2 C(\alpha_1)}{ \sphere^1}=\hmapspp{\Sigma^2 C(\alpha_1)}{ \Sigma^3 C(\alpha_1)}= \hmapspp{\Sigma^2 C(\alpha_1)}{\Sigma^5 C(\alpha_1)}= 0. \]
Similarly,
\[
\hmapspp{\Sigma^4 C(\alpha_1)}{ \Sigma^3 C(\alpha_1)}=\hmapspp{\Sigma^4 C(\alpha_1)}{\Sigma^5 C(\alpha_1)}=0.\] It follows that
\[ \hmapspp{\cp{r}{r+4}}{ \Sigma \cp{r}{r+4}} 
\cong \hmapspp{\Sigma^4 C(\alpha_1)}{\sphere^1}.\]
The group on the right-hand side above is $\bb{Z}/9$ since there is an extension
\[ 0 \rightarrow \pi_8(\bb{S}^1) \rightarrow 
\hmapspp{\Sigma^4 C(\alpha_1)}{ \sphere^1} \rightarrow \pi_4(\bb{S}^1) \rightarrow 0\]
which is nontrivial due to the Toda bracket $\langle p, \alpha_1, \alpha_1 \rangle = \alpha_2$.

When $r \equiv 3 \pmod 9$ or $r \equiv 6 \pmod 9$, the situation is slightly different. There is an additional
 $\alpha_2$ attaching map connecting the top and bottom cells of $\cp{r}{r+4}$:
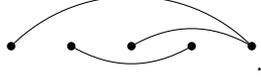
\begin{center}
\begin{tikzpicture}[scale = 0.8]
\fill (0,0) circle (2pt);
\fill (1,0) circle (2pt);
\fill (2,0) circle (2pt);
\fill (3,0) circle (2pt);
\fill (4,0) circle (2pt);
\draw (1,0) to [bend right] (3,0);
\draw (2,0) to [bend left] (4,0);
\draw (0,0) to [bend left = 45] (4,0);
\end{tikzpicture}.
\end{center}
Here, the bottom cell is again in dimension $2r$. The shorter curved line segments again denote
 $\alpha_1$-attaching maps, while the longer one denotes an $\alpha_2$-attaching map. In this case $\cp{r}{r+4}$ splits as
\[\Sigma^{2r} (Z \oplus \Sigma^2 C(\alpha_1)), \]
where  $Z$ fits into a cofiber sequence
\begin{equation}\label{eq:def_Z}\sphere^7 
\xrightarrow{\alpha_1\oplus \alpha_2} 
\sphere^4 \oplus \sphere\to Z.\end{equation}  
Thus,
\begin{align*} \hmapspp{\cp{r}{r+4}}{ \Sigma \cp{r}{r+4}} \cong & \hmapspp{\Sigma^2C(\alpha_1)}{ 
\Sigma^3C(\alpha_1)} \oplus \hmapspp{\Sigma^2C(\alpha_1)}{ \Sigma Z }\\
& \oplus  \hmapspp{Z}{ \Sigma^3C(\alpha_1)} \oplus  \hmapspp{Z}{ \Sigma Z}.\end{align*}
The calculation now requires knowledge of $\pi_i(Z)$ for $i\leq 8$. The long exact homotopy sequence associated to \eqref{eq:def_Z} shows that 
\[\pi_{ i \leq 8}(Z)\cong
\begin{cases} \Z & \text{ if } i \in\{  0, 4, 8\} \\
 \Z/3 & \text{ if } i \in\{ 3, 7\} \\
0 &\text{ otherwise.} \end{cases}\] 
A similar analysis shows that
\[ 
\hmapspp{\Sigma^2C(\alpha_1)}{\Sigma^3C(\alpha_1)} \cong  \hmapspp{\Sigma^2C(\alpha_1)}{\Sigma Z } 
\cong   \hmapspp{Z}{ \Sigma^3C(\alpha_1)} 
= 0. 
\]
We claim that $\hmapspp{Z}{ \Sigma Z} \cong \Z/3$. Indeed, 
in the exact sequence
\[ \hmapspp{\sphere^5\oplus \sphere^1}{\Sigma Z} \rightarrow \pi_7( Z)
 \rightarrow \hmapspp{Z}{ \Sigma Z} \rightarrow \hmapspp{\sphere^4\oplus \sphere}{\Sigma Z} \rightarrow 0, \]
the first homomorphism is a surjection.  
Consequently, \[\hmapspp{Z}{ \Sigma Z} \cong \hmapspp{\sphere^4\oplus \sphere}{ \Sigma Z} \cong \Z/3.\]

The remaining six cases can be similarly analyzed, so we keep the arguments brief.
\begin{itemize}
\item When $r \equiv 1 \pmod 9$, the cell structure of $\cp{r}{r+4}$ is represented as:
\begin{center}
\begin{tikzpicture}[scale = 0.8]
\fill (0,0) circle (2pt);
\fill (1,0) circle (2pt);
\fill (2,0) circle (2pt);
\fill (3,0) circle (2pt);
\fill (4,0) circle (2pt);
\draw (0,0) to [bend left] (2,0);
\draw (1,0) to [bend right] (3,0);
\end{tikzpicture}\,\,\,.
\end{center}
Therefore,
\[
\hmapsp{\cp{r}{r+4} }{\Sigma \cp{r}{r+4}} \cong 
\hmapsp{\sphere^8}{ \Sigma C(\alpha_1)}
\cong \pi_7(C(\alpha_1)) \cong \bb{Z}/9.
\]

\item When $r \equiv 4 \pmod 9$ or $r\equiv 7 \pmod 9$, the cell diagram of $\cp{r}{r+4}$ is:
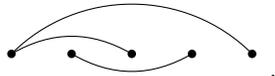
\begin{center}
\begin{tikzpicture}[scale = 0.8]
\fill (0,0) circle (2pt);
\fill (1,0) circle (2pt);
\fill (2,0) circle (2pt);
\fill (3,0) circle (2pt);
\fill (4,0) circle (2pt);
\draw (0,0) to [bend left] (2,0);
\draw (1,0) to [bend right] (3,0);
\draw (0,0) to [bend left = 45] (4,0);
\end{tikzpicture}\,\,\,.
\end{center}
Therefore, $\cp{r}{r+4}$ splits as the sum of a shift of $C(\alpha_1)$ and an indecomposable three-cell complex $Y$, so
\[\hmapspp{\cp{r}{r+4}}{\Sigma \cp{r}{r+4}} \cong \hmapspp{Y}{ \Sigma Y} \cong \Z/3.\]

\item When $r \equiv 2, 5, 8 \pmod 9$, the cell structure of $\cp{r}{r+4}$ is:
\begin{center}
\begin{tikzpicture}[scale = 0.8]
\fill (0,0) circle (2pt);
\fill (1,0) circle (2pt);
\fill (2,0) circle (2pt);
\fill (3,0) circle (2pt);
\fill (4,0) circle (2pt);
\draw (0,0) to [bend left] (2,0);
\draw (2,0) to [bend left] (4,0);
\end{tikzpicture}\,\,\,.
\end{center}
In this case, $\hmapspp{\cp{r}{r+4}}{ \Sigma \cp{r}{r+4}} \cong \pi_3(\sphere) \cong \Z/3$. 
\end{itemize}

The proof of \Cref{thm:corank_2p-2} is now complete.

\section{\texorpdfstring{$\eop$}{EO} detection for corank at most \texorpdfstring{$2p^2-p-3$}{2p²−p−3} }
\label{sec:eo_detection}

In the previous section, we computed the $p$-adic valuation of the number of isomorphism classes of stably trivial rank $r$ vector bundles on
 $\cp{r}{r+c}$ for all $r$ and $c$ such that $0<c \leq r$ and 
$c\leq 2p-3$. The next step is to study primes that are small with respect to $c$. It will be useful to have some notation: 
\begin{defn}\label{defn:phi} Let $\phrc$ denote the size of the set of isomorphism classes of rank $r$ stably trivial vector
 bundles on $\cp{}{r+c}$. Let $\phrcp$ be the largest power of $p$ dividing $\phrc$.
\end{defn}
\begin{rmk} Note that, by \Cref{subsec:background_loc}, $\phrcp$ is the order of the group
 $\hmapspp{\cp{r}{r+c}}{\Sigma \cp{r}{r+c}}$ (where all spectra are implicitly $p$-completed).\end{rmk}

As shown in \Cref{subsec:corank2p2}, the situation becomes complicated even if
$c=2p-2$. When $p$ is smaller relative to $c$, more $p$-primary attaching maps
can occur in $\cp{r}{r+c}$ and more cases must be dealt with 
because higher-degree elements in the stable homotopy groups of spheres begin to
contribute to our calculations. The idea in this section is to base change to
the category of $\eo$-modules where many such complications are suppressed since many attaching maps become trivial.

To make this precise, we need to discuss Hurewicz images. 

\begin{defn}\label{defn:Hurewicz} Given a unital ring spectrum $R$ and a spectrum $X$, the unit map $\sphere \to R$
 induces a map from stable homotopy to $R$-homology, which we denote
\[h_{R,X}\: \pi_*X \to R_*X.\]
We refer to $h_{R,X}$ as {\em the $R$-Hurewicz map for $X$}. 
\end{defn}

We are specifically interested in
\begin{equation}\label[empty]{hurewicz_surj}
\left(h_{\eo,\cp{r}{n} \otimes D\cp{r}{n} }\right)_{-1}\: \pi_{-1}(\cp{r}{n}\otimes D\cp{r}{n}) 
\to \eo_{-1}(\cp{r}{n}\otimes D\cp{r}{n}),
\end{equation}
since the order of the domain of \Cref{hurewicz_surj} is equal to $\phrcp$ and
the codomain is computable. To simplify notation, we make the following
convention.

\begin{defn}\label{def:relevant_Hurw}
When $R=\eo$ and $X=\cp{r}{n} \otimes D\cp{r}{n}$, we write $h$ for $(h_{R,X})_{-1}$ as in
 \Cref{hurewicz_surj} and will refer to $h$ simply as {\em the Hurewicz map}.
\end{defn}

The main goal of this section is to prove that $h$ is surjective if $c<2p^2-p-2$
(\Cref{thm:A}). 
The idea behind the proof of surjectivity is as follows.
In the language of \cite{Chat}, $\cp{r}{n} \otimes D\cp{r}{n}$
has {\em algebraic $\eo$-theory}, as in \Cref{def:algeo}. 
Therefore $\eo \otimes \cp{r}{n} \otimes D \cp{r}{n}$
 admits a splitting that renders
its $\eo$-theory highly computable. 
To study the image of $h$, we construct a 
cofiltration of the spectrum $\cp{r}{n} \otimes D\cp{r}{n}$ that lifts the
$\eo$-module splitting. 

Rather than studying $\cp{r}{n}\otimes D\cp{r}{n}$ in isolation,
we begin in \Cref{cofiltrations} by producing nice cofiltrations of even, finite
spectra (see \Cref{cor:exists_Xi_filtration}). 
More precisely, our cofiltrations have layers
that are even suspensions of the spectra $X_l$ defined in
\Cref{defn:alg_eo_stuff}. With this approach, we reduce surjectivity of
\Cref{hurewicz_surj} to surjectivity of 
\begin{equation}
    \label[empty]{eq:HurewiczXl2}
    h_{\eo, X_l}\: \pi_*X_l \to \eo_*X_l
\end{equation}
in certain degrees (see \Cref{lem:filter_EO_split}). 

In \Cref{computations_Xl}, we carry out the necessary computations for the
 image of \Cref{eq:HurewiczXl2} above. Combining this with results from
\Cref{cofiltrations}, we conclude that $h$ is surjective in
\Cref{prop:surjective}. 

Finally, in \Cref{lowerbound_computation}, we produce some lower bounds for
vector bundle counts.

\subsection{Preliminaries on cofiltrations}\label{cofiltrations}

We first set up a 5-lemma style argument, which will eventually reduce the
problem of surjectivity for the $\eo$-Hurewicz map for $\cp{r}{n} \otimes D\cp{r}{n}$ to a collection of simpler
problems.
\begin{lemma}
    \label{lem:filter_EO_split} 
    Let $Z$ be an even, nonzero, 
finite spectrum and $R$ a commutative ring spectrum.
Suppose that we have a finite cofiltration 
\begin{equation}
\begin{tikzcd}[row sep=1.3em]
B_0 \arrow[d,"f_0"]& B_1  \arrow[d,"f_1"] & B_2 \arrow[d,"f_2"] & & B_t \ar[d,"\simeq"]
\\ 
\mathllap{Z={}}C_0\ar[r,"g_0"]& C_1 \ar[r,"g_1"]& C_2 \ar[r,"g_2"]& \cdots \,\,  \ar[r,"g_{t-1}"]&C_t&,
\end{tikzcd}
\end{equation}
where $B_i\xrightarrow{f_i} C_i \xrightarrow{g_i} C_{i+1}$ is a cofiber sequence
for each $i$. Suppose also that:
\begin{enumerate}
\item For each $i$, the induced map on homology $H_*(f_i)$ is injective;
\item For each $i$, $f_i$ is $R$-split, meaning that $R\otimes f_i \: R\otimes B_i \to R \otimes C_i$ admits a retraction in $R$-modules; and
\item For each $i$, the map $h_{R,B_i}\:\pi_{-1}B_i \to R_{-1} B_i$ is surjective and the map $h_{R,B_i}\: \pi_{-2}B_i \to R_{-2} B_i$ is injective.
\end{enumerate}
Then $h_{R,Z}\: \pi_{-1}Z \to R_{-1}Z$ is surjective.

\end{lemma}
\begin{proof} 
By hypothesis, for the maximal cofiltration term $C_t$ we have that $C_t \simeq
B_t$. Therefore the Hurewicz map $\pi_{-1}C_t \to R_{-1}C_t$ is surjective.
Inductively, assume that we have shown the Hurewicz map $\pi_{-1}C_j\to
R_{-1}C_j$ is surjective for all $j>i$. We prove that the conclusion also holds
for $C_i$, i.e., that $h_2$ in \Cref{eq:five-lemma} below is
surjective.

For each cofiber sequence $B_i \xrightarrow{f_i} C_i \xrightarrow{g_i} C_{i+1}$,
we have long exact sequences on $\pi_*$ and $R_*$, and  a Hurewicz map between the
two long exact sequences:
\begin{equation}\label[diagram]{eq:five-lemma}
\begin{tikzcd}
\cdots \ar[r]& \pi_{-1}B_i \ar[r,"\pi_{}(f_i)"]\ar[d,"h_1"]& \pi_{-1}C_i \ar[r,"\pi_{}(g_i)"]\ar[d,"h_2"]&
 \pi_{-1}C_{i+1} \ar[r,"\delta_\pi"]\ar[d,"h_3"]& \pi_{-2} B_i\ar[d,"h_4"]\ar[r] & \cdots \\
\cdots \ar[r]&
R_{-1}B_i \ar[r,"R_{}(f_i)"]& R_{-1}C_i \ar[r,"R_{}(g_i)"]& R_{-1}C_{i+1} \ar[r,"\delta_{R}"]& R_{-2} B_i\ar[r]&\cdots &.
\end{tikzcd}
\end{equation}
Suppose that $x \in R_{-1}(C_i)$ and $R(g_i)(x)=z$. By the inductive hypothesis, we may pick a preimage
$\tilde z \in \pi_{-1}C_{i+1}$ of $z$ under $h_3$. Since $\delta_{R}(z)=0$,
$\delta_{\pi}(\tilde z)\in \ker(h_4)$. Since $h_4$ is injective, this implies
$\delta_{\pi}(\tilde z)=0$ and there exists $\tilde x \in \pi_{-1}C_i$ such that
$\pi(g_i)(\tilde x)=\tilde z$. Now, $R(g_i)(x-h_2(\tilde x)) = 0$  so
$x-h_2(\tilde x)=R(f_i)(y)$ for some $y \in R_{-1}B_i$. By hypothesis, there exists $\tilde y
\in \pi_{-1}B_i$ such that $h_1(\tilde y )=y$. Let $z:= \pi(f_i)(\tilde
y)$. Then $h_2(z) = x-h_2(\tilde x)$ so $h_2(z + \tilde x) = x$.
\end{proof}

We will apply \Cref{lem:filter_EO_split} with $Z=\cp{r}{n} \otimes D\cp{r}{n}$,
$R=\eo$.

\begin{prop}\label{cor:exists_Xi_filtration}
Let $Z$ be an even, finite spectrum. Then $Z$ admits a cofiltration 
\begin{equation}
\label[diagram]{eq:filtfilt}
\begin{tikzcd}[row sep=1.3em]
B_0 \arrow[d,"f_0"]& B_1  \arrow[d,"f_1"] & B_2 \arrow[d,"f_2"] & & B_n \ar[d,"\simeq"]
\\ 
\mathllap{Z={}}C_0\ar[r,"g_0"]& C_1 \ar[r,"g_1"]& C_2 \ar[r,"g_2"]& \cdots \,\,  \ar[r,"g_{n-1}"]&C_n&,
\end{tikzcd}\end{equation}
such that for all $i$:
\begin{itemize}
\item for some $l_i \in \{1,\ldots, p\}$ and some $s_i$ an even integer, \[B_i\simeq \Sigma^{s_i}X_{l_i},\] where
 $X_{l_i}$ is as defined in \Cref{defn:alg_eo_stuff}; and
\item $H_*(f_i)$ is split as a map of $\mathbb F_p[\xi_1]/(\xi_1^p)$-comodules.
\end{itemize}
 If additionally $Z$ has algebraic $\eo$-theory (\Cref{def:algeo}), then for
 each $i$, the map $f_i\:B_i \to C_i$ is $\eo$-split.
\end{prop}
\begin{proof}
This is a consequence of \Cref{smallcor:exists_filtration} and
\Cref{lem:hood_wants_many_lemmas} below.
\end{proof}
\begin{rmk}
Under the same hypotheses, $Z$ admits an ascending filtration with filtration
quotients suspensions of $X_l$'s.
\end{rmk}
The next lemma shows how to build a cofiltration satisfying the first item in the conclusion of \Cref{lem:filter_EO_split}.
\begin{lemma}
\label{lem:sublemma}
Let $Z$ be even, nonzero, and finite. Let $d$ be minimal such that
$H_d(Z)\neq 0$. Suppose that there is an injective map $f\:H_*(\Sigma^d X_l) \to
H_*(Z)$ of $\mathcal A_*$-comodules which is split as a map of $\mathbb
F_p[\xi_1]/(\xi_1^p)$-comodules. Then there exists a map of spectra
$\theta\:\Sigma^dX_l \to Z$ such that $H_*\theta=f$.
\end{lemma}
\begin{proof} 
Without loss of generality, assume $d=0$ so that $Z$ is connective. It suffices
to show the class of $f$ is a permanent cycle in $\op{Ext}^{0,0}(H_*(X_l),
H_*(Z))$ in the Adams spectral sequence 
\[\op{Ext}^{i+s,i}(H_*(X_l),H_*(Z))\implies \pi_{s}(\emod{\sphere}( X_l,Z)).\] 
We will show that the $E_2$ page is zero in all bidegrees which are potential
targets of differentials starting in bidegree (0,0). Since $X_l$ and $Z$ are
even spectra, there are K\"unneth and duality isomorphisms for homology so
\begin{align*}\Ext^{i+s,i}(H_*(X_l), H_*(Z))&\cong 
\Ext^{i+s,i}(\mathbb F_p, H_*(D(X_l)\otimes Z))\\
& \cong \Ext^{i+s,i}(\mathbb F_p, H_*(\Sigma^{-2(l-1)(p-1)}X_l\otimes Z)).\end{align*}
To complete the argument, we claim that 
$\Ext^{i,i+1}(\mathbb F_p, H_*(\Sigma^{-2(l-1)(p-1)}X_l \otimes Z))=0$.

Since $X_l$ is the $2(p-1)(l-1)$-skeleton of $BP$ and $Z$ is connective,
\[ 
    (X_l \otimes Z)^{2(p-1)(l-1)}=(BP \otimes Z)^{2(p-1)(l-1)}.
\]
Therefore 
\[
    \Ext^{i+s,i}\left(
        \mathbb F_p,
        H_*\left(X_l \otimes Z\right)^{2(p-1)(l-1)}
    \right)
    \cong 
    \Ext^{i+s,i}\left(
        \mathbb F_p, 
        H_*\left(BP \otimes Z\right)^{2(p-1)(l-1)}
    \right)
\]
so $ \Ext^{i+s,i}(\mathbb F_p,H_*(X_l \otimes Z))\cong \Ext^{i+s,i}(\mathbb F_p,
H_*(BP \otimes Z))$ for $s< 2(p-1)(l-1).$ Since $BP$ is complex oriented and $Z$
is even and connective, $BP\otimes Z$ splits as a sum of suspensions of $BP$.
Thus the claim follows from the fact that $\Ext^{i+s,i}(\mathbb F_p, H_*(BP))=0$
for odd $s$.
\end{proof}

\begin{cor}
    \label{smallcor:exists_filtration}
Let $Z$ be an even, finite spectrum. Then $Z$ admits a cofiltration as
in \Cref{eq:filtfilt} such that for all $i$, \[B_i\simeq \Sigma^{s_i}X_{l_i}\]
and $H_*(f_i)$ is split as a map of $\mathbb F_p[\xi_1]/(\xi_1^p)$-comodules. 
\end{cor}
\begin{proof}
The coalgebra $\mathbb F_p[\xi_1]/(\xi_1^p)$ has exactly $p$ indecomposable
comodules given by $H_*(X_l)$ for $1\leq l \leq p$. Thus $H_*(Z)$ splits as a
$\mathbb F_p[\xi_1]/(\xi_1^p)$-comodule into a sum $\bigoplus_{i} \Sigma^{s_i}
H_*(X_{l_i})$. 

Pick one of the summands with $s_i$ minimal and let
\[
f\colon \Sigma^{s_i}H_*(X_{l_i})\longrightarrow H_*(Z)
\]
denote its inclusion. Since $Z$ is even and
\[
2(l_i-1)(p-1)\leq 2(p-1)^2<2p(p-1),
\]
minimality of $s_i$ implies, for degree reasons, that the image of $f$
is an $\mathcal A_*$-subcomodule. Thus $f$ is a map of
$\mathcal A_*$-comodules. By \Cref{lem:sublemma}, there is a map
$\Sigma^{s_i}X_{l_i}\to Z$ inducing a split injection of
$\mathbb F_p[\xi_1]/(\xi_1^p)$-comodules.
The cofiber of this
map also satisfies the hypotheses of \Cref{lem:sublemma}. By induction on the
number of summands in the homology as a module over $\mathbb F_p
[\xi_1]/\xi_1^p$, we conclude the result.
\end{proof}

In order to obtain a filtration with $\eo$-split graded pieces, as in the second conclusion of \Cref{lem:filter_EO_split}, we
need the additional hypothesis that $Z$ has algebraic $\eo$-theory.

\begin{lemma}\label{lem:hood_wants_many_lemmas}
Suppose that $Z$ is an even, finite spectrum with algebraic $\eo$-theory. Let
$f\:X_l \to Z$ be a map of spectra such that $H_*(f)$ is split as a map of
$\mathbb F_p [\xi_1]/\xi_1^p$-comodules. Then $f$ is $\eo$-split: the map $\eo
\otimes f$ admits a retraction in $\eo$-modules.

\end{lemma}
\begin{proof} 
Let $\e$ denote height $p-1$ Morava $E$-theory at the ambient prime $p$.  
Given an $\eo$-module $A$, let 
\[\eeo_*(A):= \pi_*\e \otimes_{\eo}A.\]
 Note that if
$A\simeq \eo \otimes X$ for some spectrum $X$, then $\eeo_*A\cong \e_*X$
as $\eeo_*\e$-comodules.

Given $\eo$-modules $M$ and $N$, consider the relative Adams spectral sequence
for the map $\eo \to \op{E}$:
\begin{equation}
    \label[sseq]{sseq:relative} 
    E_2^{s,t}= \op{Ext}^{s,t} _{\eeo_*\e}(\eeo_*M, \eeo_*N) \implies \pi_{t-s}\eomod(M,N).
\end{equation}
By \cite[Corollary 3.4]{Devinatz}, if $M$ is finite and relatively projective,
then the relative Adams spectral sequence is strongly convergent.

The $\eo$-module $\eo \otimes Z$ is finite since $Z$ is a finite spectrum. Since
$Z$ is even and $\e$ complex orientable, $\e \otimes Z$ is a direct sum of
shifts of $\e$ and $\e_*Z$ is projective over $\e_*$. Specializing to 
$M=\eo \otimes Z$ and $N=\eo \otimes X_l$ in \Cref{sseq:relative} gives a strongly
convergent spectral sequence
\begin{equation}
    \label[sseq]{sseq:ZXl} 
    E_2^{s,t}= \op{Ext}^{s,t} _{\eeo_*\e}(\e_*Z, \e_*X_l) \implies \pi_{t-s}\eomod(\eo\otimes Z,\eo\otimes X_l).
\end{equation}

We now produce a retraction $g\: \eo \otimes Z \to \eo \otimes X_l$ of $\eo
\otimes f\: \eo \otimes X_l \to \eo \otimes Z$. We first prove the case $l=p$ and then prove $l<p$ by downwards
induction on $l$. In fact, we prove by downward induction the stronger assertion that the retraction can be chosen to lift the projection onto the specified $H_*(X_l)$-summand in the chosen homological decomposition. So, we fix a decomposition $\eo \otimes Z \cong \oplus_i\Sigma^{s_i}\eo \otimes X_{l_i}$ in $\eo$-modules lifting a given homology decomposition.

Suppose $f\: X_p \to Z$ is a given map inducing an injection on
homology.  By \cite[Proposition 5.20]{Chat} the edge morphism in \Cref{sseq:ZXl}, i.e., the map
\[\pi_0\eomod(\eo \otimes Z, \eo \otimes X_p) \to \op{Hom}_{\eeo_*\e}(\e_*(Z), \e_*(X_p))\]
given by $p \mapsto \pi_*(\e \otimes_{\eo} p)$, is an isomorphism. Therefore, it suffices to produce a splitting 
\[\bar g \in \op{Hom}_{\eeo_*\e}(\e_*(Z), \e_*(X_p))\] of $\pi_*(\e \otimes_{\eo}f).$
The existence of $\bar g$ follows from the fact that $\e_*(X_p)$ is a relatively injective $\eeo_*\e$-comodule. Moreover, we choose this splitting to be compatible with the chosen decomposition $\eo \otimes Z \cong \oplus_i\Sigma^{s_i}\eo \otimes X_{l_i}$ which exists and lifts the $\eo$-homology decomposition of $Z$, since $Z$ has algebraic $\eo$-theory (see \cite[Theorem 5.6]{Chat}).

Now assume $l<p$. We are given $f\: X_l \to Z$ with $H_*(f)$ split as a map of
$\mathbb F_p[\xi_1]/\xi_1^p$-comodules. Let $N:=2(p-1)(l-1)+2p-3$. There is a
generator $\alpha_1^{(l)}\: \sphere^N \to X_l$ such that $\cof(\alpha_1^{(l)})\simeq
X_{l+1}$. Let $Z':=\cof(f\circ \alpha_1^{(l)})$ so that we have $f'\:X_{l+1}\to Z'$
and a commutative diagram: 

\begin{equation}\label[diagram]{inductivestep}
\begin{tikzcd}
\sphere^N \ar[d]\ar[r,"="] & \sphere^N\ar[d]\\
X_l \ar[r, "f"]\ar[d] & Z\ar[d] \\
X_{l+1} \ar[r, "f'"] & Z'\,,
\end{tikzcd}
\end{equation}
where the vertical sequences are both cofiber sequences. By \cite[Lemma 5.11]{Chat}, $Z'$ has algebraic
 $\eo$-theory. Note also that $f'$ induces a split map of $\mathbb F_p[\xi_1]/\xi_1^p$-comodules.
 By the inductive hypothesis, we have a splitting of $\eo$-modules
\[ \begin{tikzcd} \eo \otimes X_{l+1} \ar[r,"\eo \otimes f'"] & \eo \otimes Z'\ar[l,"g' ", bend left]
\end{tikzcd}
\]
which is compatible with the identification $\eo \otimes Z'\simeq \oplus_i \Sigma^{s_i}(\eo \otimes X_{l_i})$ arising from a chosen $\eo$-homology splitting of $Z'$ (see \cite[Theorem 5.6]{Chat}).

Now rotate \Cref{inductivestep} and tensor with $\eo$ to obtain a commutative diagram:

\begin{center}
\begin{tikzcd}
\eo\otimes  X_l \ar[r, "f"]\ar[d] & \eo\otimes Z\ar[d] \\
\eo\otimes X_{l+1}\ar[d,"p_1"] \ar[r, "f'"] & \eo\otimes Z'\ar[d,"p_2"]\\
\eo\otimes \sphere^{N+1}\ar[r,"="] & \eo\otimes \sphere^{N+1} & .
\end{tikzcd}
\end{center}
The splitting of $\eo \otimes f'$ together with the fact that the splitting identifies $\eo \otimes X_{l+1}$ with a summand in the $\eo$-module decomposition of $\eo \otimes Z'$ gives an identification $\eo \otimes Z'
\simeq M \oplus (\eo \otimes X_{l+1})$, and $f'$ is identified with the inclusion of the second summand. The arrow labeled $p_2$ above is
projection onto $\eo\otimes X_{l+1}$ composed with the map $p_1\: \eo\otimes
X_{l+1} \to \eo \otimes S^{N+1}$. Therefore $\fib(p_2) \simeq M \oplus (\eo
\otimes X_l)$ and by commutativity $f$ is the inclusion of the summand $\eo\otimes
X_l$ as was to be shown.
\end{proof}

\begin{rmk}
    \label{rmk:homology_support} 
    Note that $X_l$ has homology concentrated in degrees $0,\ldots, 2(l-1)(p-1)$.
    Therefore if $H_*(Z)$ is concentrated in an interval $[n,m]$, each  $s_i$ in
    \Cref{cor:exists_Xi_filtration} satisfies $n\leq s_i\leq m-2(l_i-1)(p-1)$.
\end{rmk}

\subsection{Computations for suspensions of \texorpdfstring{$X_l$}{X\_l}}\label{computations_Xl}
In this section, we apply \Cref{lem:filter_EO_split} in the case of interest. Let $Z= \cp{r}{n}\otimes D\cp{r}{n}.$
 \Cref{cor:exists_Xi_filtration} implies that $ Z$ is cofiltered by spectra $\Sigma^{2s}X_l$ such that
 each cofiltration quotient is $\eo$-split. 
Thus, by \Cref{lem:filter_EO_split}, we can study surjectivity of the Hurewicz map $\pi_{-1}Z \to \eo_{-1}Z$ 
by understanding the Hurewicz maps 
$\pi_*\Sigma^{2s} X_l \to {\eo}_* \Sigma^{2s}X_l.$ 

First, we prove the following lemma:
\begin{lemma}\label{lem:map_sseq} Let $ A_{*,*}^n $ and $ B_{*,*}^n$ be Adams-graded homological spectral sequences.
 Let \[\{f_{*,*}^n\: A_{*,*}^n \to B_{*,*}^n\}\] be a map of spectral sequences. Suppose that there is an odd integer
 $k$ and an even integer $k'>k$ such that, for all integers $t$ such that $k< t < k'$, $f_{t,*}^2$ is injective for $t$
 even and surjective for $t$ odd. 
Then, for all integers $n$ and $t$ such that $n \geq 2$ and $k<t<k'$, $f_{t,*}^n$ is injective for
 each even $t$ and is surjective for each odd $t$. \end{lemma}
\begin{proof} A differential on the $r$-th page with source $(t,s)$ has target $(t-1,s-r)$. 
Suppose, for a contradiction, that $n\geq 2$ is minimal such that $f_{t,*}^{n+1}$ 
fails to be injective for even $t$ or fails to be surjective for odd $t$.

Suppose the map $f_{t,s}^{n+1}\: A_{t,s}^{n+1} \to B_{t,s}^{n+1} $ is
non-injective for some even $t$ with $k < t < k'$ and some $s \in \mathbb Z$.
Pick a nonzero $x \in A_{t,s}^{n+1}$ such that $f^{n+1}_{t,s}(x)=0$. Let $\tilde{x}$ be a representative for $x$ on the $n$-th page. Then $f^n_{t,s}(\tilde{x})$
is the target of a differential with source $y\in B_{t+1,s+n}^{n}$.  
Since $k'$ is even and $k<t<k'$, we have that $k<t+1<k'$. 
By surjectivity on the $n$-th page, there exists $y'\in A_{t+1,s+n}^{n}$ with $f^n_{t+1,s+n}(y')=y$. 
By naturality, $d(y')-\tilde{x}$ is in the kernel of $f_{t,s}^{n}$. By injectivity at the $n$-th stage, 
$d(y')=\tilde{x}$ and $\tilde{x}$ represents the zero class on the $E_{n+1}$-page. Hence $x=0$ on the $(n+1)$-page, a contradiction.

If the map $f_{t,s}^{n+1}$ fails to be surjective for some odd $t$ with $k<t<k'$ 
and some $s \in \mathbb Z$, then there is a class $y\in B_{t,s}^{n+1}$ that is not in the image of 
$f_{t,s}^{n+1}$. Since $f_{t,s}^n$ is surjective, there is a class $x \in A_{t,s}^n$ with $f_{t,s}^n(x)=\tilde y$, where $\tilde y$ is a representative for $y$ on the $n$-th page. 
For surjectivity to fail on the $n+1$-page, the class $x$ must support a nonzero differential $d_{t,s}^n\: A_{t,s}^n \to A_{t-1,s-n}.$ 
By naturality and the fact that  $\tilde y$ survives to the $n+1$-page:
\begin{align*}
f^n_{t-1,s-n}(d^n_{t,s}(x))&=d^n_{t,s}(\tilde{y})=0. \end{align*}
Since $k<t<k'$, $k$ is odd, and $t$ is odd, we have that $t-1$ is even and $k< t-1 < k'$. Therefore,
$d^n_{t,s}(x)=0$ by injectivity of $f^n_{t-1, s-n}$, so the class $y$ is in
the image of $f^{n+1}$.
\end{proof}

\begin{cor}\label{cor:inj_surj_xl}  Let $A=\Sigma^{2s}X_{l}$. If $p$ is an odd prime and 
\[-\left(2p^2-p-2\right)< s < (p^2-1)(p-1)-(l-1)(p-1)-2,\] the Hurewicz map $\pi_{-1}A\to \eo_{-1}A$ 
is surjective and the Hurewicz map $\pi_{-2}A \to \eo_{-2}A$ is injective. 
\end{cor}
\begin{proof}
Consider the homological Atiyah--Hirzebruch spectral sequences 
\begin{equation}\label[sseq]{AHSS1}A_{*,*}^2=H_q(X_l; \pi_{t-q}\sphere ) \implies \pi_{t}(X_l)\end{equation}
and
\begin{equation}\label[sseq]{AHSS2}B_{*,*}^2=H_q(X_l; \pi_{t-q}\eo)\implies \eo_t(X_l).\end{equation}
\Cref{AHSS1} and \Cref{AHSS2} are related by a Hurewicz map $h_{*,*}\: A_{*,*}^n \to B_{*,*}^n$, given on the
$E_2$-page by 
\begin{equation}\label{eq:e2page}y \otimes x \mapsto h_{\eo, \sphere}(y)\otimes x\end{equation} 
for $x$ a homology class and
$y \in \pi_*\sphere$. 
We identify $k,k'$ such that $h$ as in \eqref{eq:e2page} satisfies the hypotheses of \Cref{lem:map_sseq}.
We claim that

\begin{equation}\label[empty]{Hurewicz_pt}h_{\eo,\sphere}\: \pi_{t-q}\sphere \to \eo_{t-q}\sphere\end{equation}
is injective on even degrees for all
$t-q<4p^2-2p-4.$ 
The first two even-degree $p$-torsion elements in
$\pi_*\sphere$ are $\beta_1$ and $\beta_2$ \cite[4.4.20]{Ravenel_green}. 
The smallest positive even-degree element in the kernel of \Cref{Hurewicz_pt} is $\beta_2$ in degree $4p^2-2p-4$. By
\cite[Theorem 2.1]{Nave}, all odd-degree homotopy classes in $\pi_*\eo$ for 
$0
\leq * \leq |\beta_2|  $ are in the image of \Cref{Hurewicz_pt}.

To choose the lower bound $k$ for an odd prime, we note that the largest $i < 0$ such that 
$\pi_i\eo\neq 0$ is $i = -2(p^2-1)(p-1)$. The group $\pi_{-2(p^2-1)(p-1)}$ is generated by 
$\Delta^{-p}\cdot \beta_1^{(p-1)^2}$.

In summary, for odd primes $p$, the Hurewicz map $h_{\eo,\sphere}\:\pi_*\sphere \to \eo_*\sphere$ is
injective on even degrees and surjective on odd degrees for
$-2(p^2-1)(p-1)<*<4p^2-2p-4$.

Thus, 
\[h_{*,*}^2\:\pi_{t-q}\sphere\otimes H_q(X_l) \to \pi_{t-q}\eo \otimes
H_q(X_l)\] is injective for even $t-q$ and surjective for odd $t-q$ as long as
\[-2(p^2-1)(p-1)<t-q< 4p^2-2p-4.\]
The spectrum $X_l$ has cells in dimensions $(2p-2)i$ where $0\leq i \leq l-1$.
The nonzero terms on the $E_2$-page of both spectral sequences
lie in bidegrees $(t,q)$ with $q=(2p-2)i$ for $0 \leq i \leq l-1$. 
Take
\begin{align*}
k&=-2(p^2-1)(p-1)+2(l-1)(p-1)+3\\
k'&=4p^2-2p-4\end{align*} in \Cref{lem:map_sseq}. We have shown that
$h_{t,*}\: A_{t,*}^2 \to B_{t,*}^2$ is injective on even $t$ and surjective
on odd $t$ for $k<t<k'.$ 
By \Cref{lem:map_sseq}, the Hurewicz map $h_{\eo, X_l}$ is injective in degrees $-2-2s$ for $k<-2-2s<k'$, and
 surjective in degrees $-1-2s$ for $k< -1-2s < k'.$
\end{proof}

\begin{prop}\label{prop:surjective}
    For $p$ odd and $c< 2p^2-p-2,$ the map
    \[\pi_0\smod(\cp{r}{n},\Sigma\cp{r}{n}) \to \pi_0\eomod(\eo\otimes \cp{r}{n}, \eo \otimes \Sigma\cp{r}{n})\]
    is surjective.\end{prop}
\begin{proof}
    \Cref{cor:exists_Xi_filtration} and \Cref{rmk:homology_support} imply that
    $\cp{r}{n}\otimes D\cp{r}{n}$ is cofiltered by spectra $\Sigma^{2s}X_l$ where
    $1\leq l \leq p$ and $-c\leq s \leq c-(l-1)(p-1)$ and that each cofiltration
    quotient is $\eo$-split. 
   For $p$ odd, if $c < 2p^2-p-2$ then, by \Cref{cor:inj_surj_xl}, the hypotheses of \Cref{lem:filter_EO_split} are satisfied. 
\end{proof} 
\begin{rmk}
    For odd $p$, the surjectivity range for the Hurewicz image $\pi_{-1}\Sigma^{2s}X_l \to
    \eo_{-1}\Sigma^{2s}X_l$ given in \Cref{cor:inj_surj_xl} can be extended to
    have lower bound $-\frac{|\alpha_1\Delta|}{2}+1<s$ using an Adams spectral
    sequence argument. However, this does not improve our results on vector
    bundle counts, since injectivity is also needed to reduce the Hurewicz image
    for $\cp{r}{n} \otimes D \cp{r}{n}$ to the Hurewicz image for even
    suspensions of $X_l$'s. 
\end{rmk}
\begin{rmk}\label{remark:bijection} The methods used to prove \Cref{prop:surjective} above can be modified to show that the Hurewicz map in question is bijective for $n-r\leq 2p-3$. However, we have already established the corresponding enumeration result via direct computation (cf. \Cref{thm:count1}), so we do not give the argument. 
\end{rmk}

\subsection{Lower bounds for vector bundle counts and \texorpdfstring{$\eo_*X_l$}{EO*X\_l}}\label{lowerbound_computation}

\Cref{prop:surjective} implies that we can give a lower bound for $\phi(r,n-r,p)$ when 
\begin{equation}\label{eq:relevant_corank} 2p-3< n-r < \frac{2p^2-p-2}{2},\end{equation} 
by computing $\eo_{-1} \cp{r}{n} \otimes D\cp{r}{n}$ in full. For the range of coranks as in 
\Cref{eq:relevant_corank}, we can compute the relevant 
$\eo$-homology in full. The starting point is that $\cp{r}{n}$ has algebraic 
$\eo$-theory \cite{Chat}, so $\eo\otimes \cp{r}{n}\otimes D\cp{r}{n}$ 
splits as a sum of even suspensions of $\eo \otimes X_l$ for $1 \leq l \leq p$. 
The exact form of this splitting is determined by combinatorial formulas in \cite{Chat}.
Therefore, we will build up a complete computation of the $\eo$-theory of $\cp{r}{n}\otimes D\cp{r}{n}$ from the following result.
\begin{prop}\label{prop:alphal} For $p$ an odd prime and $-\frac{(2p^2-p-2)}{2}< s < \frac{2p^2-p-2}{2}-(p-1)(l-1)$, $1\leq l \leq p$,
\[\eo_{-1}\Sigma^{2s}X_l \cong \begin{cases} \Z/p \text{ if $l \neq p$ and $s=-l(p-1)$}\\ 0 \text{ otherwise.}\end{cases}\]
\end{prop}
\begin{proof}
If $l=p$, then by \cite[Corollary 3.5]{Chat},
\[
\eo\otimes X_p\simeq \e^{hC_{(p-1)^2}}.
\]
Since $(p-1)^2$ is prime to $p$ and $\pi_*\e$ is concentrated in even
degrees, the homotopy fixed-point spectral sequence shows that
$\pi_*(\eo\otimes X_p)$ is concentrated in even degrees. Thus
$\eo_{-1-2s}X_p=0$ for every integer $s$. We may therefore assume that
$l<p$.

By \Cref{cor:inj_surj_xl}, the Hurewicz map 
$\pi_{-1-2s}X_l \to \eo_{-1-2s}X_l$ is surjective. By \cite[Proposition 5.9]{Chat}, 
$\pi_{-1+2l(p-1)}\eo \otimes X_l$ is generated by the 
$p$-torsion elements $\alpha_1^{(l)}$, where $\alpha_1^{(l)}$ is the Hurewicz image of the attaching map for building the complex $X_{l+1}$. Moreover, $\pi_{-1+2k(p-1)}\eo \otimes X_l =0$ for 
$k\neq l$. This completes the computation for $s$ a multiple of $(p-1)$.

To see that all other groups $\pi_{-1-2s}\eo \otimes X_l$ are zero, consider 
the Atiyah--Hirzebruch spectral sequence $H_q(X_l;\pi_{t-q}\eo)\implies \eo_tX_l$.
The only possible nonzero terms on the $E_2$-page are in groups 
\[H_{2i(p-1)}(X_l;\pi_{t-2i(p-1)}\eo)\] for $q=2i(p-1)$ with $0 \leq i \leq l-1$.  
For $t$ odd and \[-(2p^2-p-2)\leq t-2i(p-1)< 2p^2-p-2,\] the only 
possible classes are given on the $E_2$-page by $\alpha_1$ on a cell in $X_l$, which are on the 
vertical line $t=-1+2k(p-1)$ for some $k$.
\end{proof}
We next compute $\eo_{-1}(\cp{r}{n} \otimes D\cp{r}{n})$ in the cases of interest. 
However, to state the result, we need a few definitions.
\begin{defn} Given $k \in \mathbb Z$ and $d$ a positive integer, let 
$[k]_d=k \pmod d\in \{0,\ldots, d-1\}.$ For integers $n$ and $j$, let $n_j =n-[n-j]_{p-1}$. 
Explicitly, $n_j$ is the largest integer less than or equal to $n$ with the same remainder as $j$ $\pmod{p-1}$.
\end{defn}

\begin{thm}\label{thm:total} Let $p$ be an odd prime and let $n-r < \frac{2p^2-p-2}{2}.$ Let $j(n,r)$ be the number of $i$ in $\{0, \ldots, p-2\}$ such that all of the following hold:
\begin{itemize}
\item[{\bf A}.]  $ r+i + ([r+i]_p)(p-1) < n_{r+i} ;$ 
\item[{\bf B}.] $[r+i]_p+1\leq [n_{r+i}]_p$; and
\item[{\bf C}.] $n_{r+i}-(r+i)+(p-1)\left([n_{r+i}]_p-[r+i]_p\right) = p(p-1).$
\end{itemize}
Then $\eo_{-1}(\cp{r}{n} \otimes D\cp{r}{n})\cong \left(\Z/p\right)^{\oplus j(n,r)}.$ 
\end{thm} \begin{rmk} In the larger range $n-r<2p^2-p-2$, the proof of \Cref{thm:total} in fact shows that $\eo_{-1}(\cp{r}{n} \otimes D\cp{r}{n})$ contains a summand isomorphic to $(\Z/p)^{\oplus j(n,r)}$, with the same definition of $j(n,r)$. Hence, we can still deduce a lower bound for $\pi_{-1}\cp{r}{n}\otimes D\cp{r}{n}$ by \Cref{prop:surjective}. The smaller range for the full computation of $\eo_{-1}(\cp{r}{n} \otimes D\cp{r}{n})$ is a consequence of the restrictions introduced in the Atiyah--Hirzebruch spectral sequence argument for giving a complete computation of $\eo_{-1}\Sigma^{2s}X_l$ in \Cref{prop:alphal}; a partial computation of  $\eo_{-1}\Sigma^{2s}X_l$ that still picks out a large summand in $\eo_{-1}(\cp{r}{n} \otimes D\cp{r}{n})$ follows directly from \cite[Proposition 5.9]{Chat}.\end{rmk}

We refer the reader to \Cref{combinatorics} for a proof of \Cref{thm:total}, and turn our attention here to applications. First, we note that, while \Cref{prop:surjective} guarantees that $\pi_0\smod(\cp{r}{n},\Sigma\cp{r}{n})$ surjects onto $\eo_{-1}\cp{r}{n} \otimes D\cp{r}{n}$ for $n-r<2p^2-p-2$ and hence also for $n-r<\frac{2p^2-p-2}{2}$, condition (C) in \Cref{thm:total} above will frequently result in $\eo_{-1}\cp{r}{n} \otimes D\cp{r}{n}$ being zero, meaning that we do not actually get an interesting lower bound for the number of stably trivial non-trivial rank $r$ vector bundles on $\cp{}{n}$. However, we work this out for $p=3$, and demonstrate that some new information is obtained. We then show that the results are more striking for $p$ large, where there is more flexibility to satisfy (A) - (C).

First, consider $j(n,r)$ for $p=3$.
\begin{cor}\label{cor:p3-case} If $p=3$ and $r\geq 3$, the values of $j(n,r)$ are as indicated below:
 \vspace{.5em}
    \begin{center}
        \begin{tabular}{|c|c|c|c|c|c|c|c|}
        \hline
            $n-r$ & $2$ & $3$ & $4$ & $5$ & $6$ \\
        \hline
            $j(n,r)$ for $n\equiv 0\pmod{3}$  & $0$ & $1$ & $1$ & $1$
& $1$ \\
       $j(n,r)$ for $n\equiv 1\pmod{3}$  & $0$ & $0$ & $1$ & $1$ 
& $0$  \\
       $j(n,r)$ for $n\equiv 2\pmod{3}$  & $1$ & $1$ & $1$ & $2$ 
& $1$  \\
        \hline
        \end{tabular}.
    \end{center}
  \vspace{.5em}
All values are zero for $n-r>6$.
 \end{cor} 

\begin{proof} First, note that $n_k=n$ if $k \equiv n \pmod{2}$, and $n_k=n-1$ if $k\not\equiv n \pmod{2}$.
\begin{itemize}
\item[(1a)] Suppose $i=0$, $r\equiv n\pmod{2}$. Then (A)-(C) become:
\begin{itemize}
\item[(A1a)] $r+2[r]_3<n$;
\item[(B1a)] $[r]_3+1\leq [n]_3$; and
\item[(C1a)] $n-r+2([n]_3-[r]_3)=6$.
\end{itemize}
The inequality (B1a) has a solution when either $n\equiv 1\pmod{3}$ and $r\equiv 0\pmod{3}$, or $n\equiv 2\pmod{3}$ and $r\equiv 0\text{ or }1 \pmod{3}$.
For these values of $[n]_3$ and $[r]_3$, we find that (C1a) has a solution precisely if:
\begin{itemize}
\item  $n\equiv 1 \pmod{3}$, $n-r=4$ (which forces $r\equiv 0\pmod{3}$); or
\item $n\equiv 2 \pmod{3}$, $n-r=2$ (which forces $r\equiv 0 \pmod{3}$); 
\item $n\equiv 2 \pmod{3}$, $n-r=4$ ($r\equiv 1 \pmod{3}$).
\end{itemize}
In all of these cases (A1a) is satisfied.
\item[(1b)] Suppose $i=0$, $r\not\equiv n\pmod{2}$. Then (A)-(C) become:
\begin{itemize}
\item[(A1b)] $r+2[r]_3<n-1$;
\item[(B1b)] $[r]_3+1\leq [n-1]_3$; and
\item[(C1b)] $n-1-r+2([n-1]_3-[r]_3)=6$.
\end{itemize}
The inequality (B1b) has a solution when either $n-1\equiv 1\pmod{3}$ and $r\equiv 0\pmod{3}$, or $n-1\equiv 2\pmod{3}$ and $r\equiv 0\text{ or }1 \pmod{3}$.
For these values of $[n-1]_3$ and $[r]_3$, we find that (C1b) has a solution precisely if:
\begin{itemize}
\item  $n\equiv 2 \pmod{3}$, $n-r=5$ (which forces $r\equiv 0\pmod{3}$); or
\item $n\equiv 0 \pmod{3}$, $n-r=3$ (which forces $r\equiv 0 \pmod{3}$); or
\item $n\equiv 0 \pmod{3}$, $n-r=5$ (which forces $r\equiv 1 \pmod{3}$).
\end{itemize}
In all of these cases, (A1b) is satisfied.
\item[(2a)] Suppose $i=1$, $r\equiv n\pmod{2}$. Then $n_{r+1}=n-1$. Then (A)-(C) become:
\begin{itemize}
\item[(A2a)] $r+1+2[r+1]_3<n-1$;
\item[(B2a)] $[r+1]_3+1\leq [n-1]_3$; and
\item[(C2a)] $n-1-(r+1)+2([n-1]_3-[r+1]_3)=6$.
\end{itemize}
The inequality (B2a) has a solution when either $n-1\equiv 1\pmod{3}$ and $r+1\equiv 0\pmod{3}$, or $n-1\equiv 2\pmod{3}$ and $r+1\equiv 0\text{ or }1 \pmod{3}$.
For these values of $[n-1]_3$ and $[r+1]_3$, we find that (C2a) has a solution precisely if:
\begin{itemize}
\item  $n\equiv 0 \pmod{3}$, $n-r=4$ (which forces $r\equiv 2\pmod{3}$); or
\item $n\equiv 0 \pmod{3}$, $n-r=6$ (which forces $r\equiv 0 \pmod{3}$); 
\item $n\equiv 2 \pmod{3}$, $n-r=6$ ($r\equiv 2 \pmod{3}$).
\end{itemize}
In all of these cases (A2a) is satisfied.
\item[(2b)] Suppose $i=1$, $r\not\equiv n\pmod{2}$. Then $n_{r+1}=n$. Then (A)-(C) become:
\begin{itemize}
\item[(A2b)] $r+1+2[r+1]_3<n$;
\item[(B2b)] $[r+1]_3+1\leq [n]_3$; and
\item[(C2b)] $n-(r+1)+2([n]_3-[r+1]_3)=6$.
\end{itemize}
The inequality (B2b) has a solution when either $n\equiv 1\pmod{3}$ and $r+1\equiv 0\pmod{3}$, or $n\equiv 2\pmod{3}$ and $r+1\equiv 0\text{ or }1 \pmod{3}$.
For these values of $[n]_3$ and $[r+1]_3$, we find that (C2b) has a solution precisely if:
\begin{itemize}
\item  $n\equiv 1 \pmod{3}$, $n-r=5$ (which forces $r\equiv 2\pmod{3}$); or
\item $n\equiv 2 \pmod{3}$, $n-r=3$ (which forces $r\equiv 2 \pmod{3}$); 
\item $n\equiv 2 \pmod{3}$, $n-r=5$ ($r\equiv 0 \pmod{3}$).
\end{itemize}
In all of these cases (A2b) is satisfied.

\end{itemize}
Summarizing the lists from each case, we get the claimed result.
\end{proof} 

\begin{rmk} Let $p=3$. At corank $2$ and $3$, $\eo$-theory $\eo_{-1}\cp{r}{n} \otimes D\cp{r}{n}$ is exactly the $3$-component of $\pi_0\smod(\cp{r}{n},\Sigma\cp{r}{n})$, as noted in \Cref{remark:bijection}. In corank $2$, \Cref{cor:p3-case} recovers the $3$-divisibility pattern in \cite[Theorem 1.2]{Hu}. In corank $3$, compare with \Cref{thm:count1} and \Cref{cor:bundle_ct1}: taking $n=r+c$, the condition for a non-zero exponent to appear is that one of $n-2$ and $n-3$ is divisible by $3$, which matches the answer in the table above. \end{rmk} 

For arbitrary odd $p$ and appropriately chosen $c$, the results above are more surprising. We illustrate this by computing $\eo_{-1}\cp{r}{n}\otimes D\cp{r}{n}$ for $n-r=(p-1)^2$. 
\begin{cor}\label{cor:bigcount} For each odd prime $p$, each integer $t\geq 0$, and $0 \leq r \leq p-1$, 
the number of rank $(p^2+tp+r)$-bundles on $\cp{}{(p-1)^2+ p^2+tp+r}$ with vanishing Chern classes is a non-zero multiple of $p^{K}$, where
\begin{equation}\label{eq:first3}
K=\begin{cases} p-2& \text{ if } [r-1]_p\geq p-3 \\
p-3  &  \text{ if } [r-1]_p< p-3.
\end{cases}
\end{equation}
\end{cor}
\begin{proof} Note that the corank is fixed at $(p-1)^2< \frac{2p^2-p-2}{2}$. For any $t\geq 0$ and $0 \leq r \leq p-1$, we are in the metastable range.
By \Cref{prop:surjective}, the order of the group 
\[\eo_{-1}\left(\cp{p^2+pt+r}{(p-1)^2+p^2+pt+r} \otimes D\cp{p^2+pt+r}{(p-1)^2+p^2+pt+r}\right)\]
 is a divisor of the number of stably trivial rank $p^2+tp+r$ bundles on $\cp{}{p^2+(p-1)^2+tp+r}$. 

By \Cref{cor:eo_periodicity}, 
\[\eo\otimes\left( \cp{p^2+pt+r}{(p-1)^2+p^2+pt+r}\otimes
D\cp{p^2+pt+r}{(p-1)^2+p^2+pt+r}\right)
\simeq \eo\otimes\left( \cp{r}{(p-1)^2+r}\otimes
D\cp{r}{(p-1)^2+r}\right).\]
We apply \Cref{thm:total} with $n=(p-1)^2+r$. Since $n \equiv r \pmod{p-1}$, we
obtain:
\begin{align*} 
n_{r}& =(p-1)^2+r
& [n_{r}]_p& = [r+1]_p
\\
n_{r+1}&=(p-1)^2-(p-1)+(r+1)&
[n_{r+1}]_p&=[r+3]_p\\
n_{r+2}&=(p-1)^2-(p-1)+(r+2)
&[n_{r+2}]_p&=[r+4]_p
\\
&\vdots  & \vdots \\
 n_{r+p-2}&=(p-1)^2-(p-1)+(r+p-2) 
& [n_{r+p-2}]_p&=[r+p]_p.
\end{align*}

That is, for $0 \leq r \leq p-1$ and $0\leq i \leq p-2$,
\begin{align}\label{eq:nri}
n_{r+i}&=\begin{cases} (p-1)^2+r & \text{ if } i=0\\
(p-1)^2-(p-1)+(r+i)& \text{ if } i\neq 0.
\end{cases}\\ \nonumber
[n_{r+i}]_p&=\begin{cases} [r+1]_p &  \text{ if } i=0\\
[r+i+2]_p &\text{ if } i\neq 0.
\end{cases}
\end{align}
We first deal with $i>0$. In this case, conditions (A) - (C) of \Cref{thm:total} become
\begin{enumerate}
\item[{ A}.]  $ ([r+i]_p)(p-1) < (p-1)^2-(p-1)=(p-2)(p-1)$; 
\item[{ B}.] $[r+i]_p+1\leq [r+i+2]_p$ ; and
\item[{ C}.] $(p-1)^2-(p-1)+(p-1)\left([r+i+2]_p-[r+i]_p \right) = p(p-1).$
\end{enumerate}
Condition (A) is satisfied if and only if $0\leq [r+i]_p <p-2.$ 
In this range, $[r+i+2]_p=[r+i]_p+2$, so condition (B) is satisfied.

Now consider condition (C). Substituting $[r+i+2]_p=[r+i]_p+2$ shows:
\begin{align*}(p-1)^2-(p-1)+(p-1)\left([r+i]_p+2-[r+i]_p \right) &= (p-1)^2-(p-1)+2(p-1)\\
&= (p-1)(p-1-1+2)\\
&=p(p-1),\end{align*}
so (C) is satisfied.

Thus, for each $1\leq i\leq p-2$ such that $[r+i]_p<p-2,$ we get a summand isomorphic to $\Z/p$ in the group of interest. If $r\equiv -1\pmod{p}$, there are $p-2$ such $i$. If $r\equiv 0$ or $r \equiv -2 \pmod{p}$, there are $p-3$ such $i$. Otherwise, there are $p-4$ such $i$. 

We next consider the case $i=0$. In this case, conditions (A) - (C) become:
\begin{itemize}
\item[{A}.]  $ ([r]_p)(p-1) < (p-1)^2$;
\item[{B}.] $ [r]_p+1\leq [r+1]_p$; and
\item[{C}.] $(p-1)^2+(p-1)\left([r+1]_p-[r]_p\right) = p(p-1).$
\end{itemize}
Item (A) is equivalent to 
$ 0 \leq [r]_p < p-1.$
If this is satisfied, $[r+1]_p=[r]_p+1$, and so condition (B) is satisfied. Equation (C) is automatically satisfied.
So, we get a $\Z/p$ summand precisely when $0\leq [r]_p<p-1$. 

Thus, we get one additional summand unless $r \equiv -1\pmod{p}$. 

In total, this shows that we have $p-2$ summands when $r\equiv -2,-1,0\pmod{p}$, and $p-3$ summands otherwise. We can rephrase this as $p-2$ summands when $[r-1]_p\geq p-3$, and $p-3$ otherwise. 
\end{proof}
The reader will note that the methods used to deduce 
\Cref{cor:bigcount} above could be applied to obtain many other families of examples from \Cref{thm:total}.
As a special case of \Cref{cor:bigcount}, we obtain \Cref{thm:X}.

\subsection{Geometric interpretations of vector bundles detected in $\eo$-theory at corank $2p-2\leq c < 2p^2-p-2$}\label{geometric-2}

We now note that a story analogous to that explained in \Cref{geometric-easy} can be partially adapted to higher corank $2p-2 \leq c < 2p^2-p-2.$ However, the proofs of these observations do require a somewhat detailed study of the proof of \Cref{thm:total}. 

Given \Cref{geometric-easy}, a natural question is whether the vector bundles detected in higher corank can be constructed as pullbacks of vector bundles on spheres. We show that this is sometimes, but rarely, the case. We formulate analogues of the bundles defined by \Cref{diagram-V}.

Suppose that $\cp{r}{n}$ has $(r+i)$-th Adams summand $A_{i}$ with top cell in degree $n_{r+i} \equiv p-1 \pmod{p}$. Then $\eo \otimes A_i$ has a minimal-degree $\eo$-split summand $\Sigma^{2(r+i)}\eo \otimes X_{l(i)}$. For $l(i)<p$, we have a class 
\[\alpha_1^{(l(i))}\:\sphere^{(l(i)-1)(2p-2)+2p-3}\to X_{l(i)}\] 
which is the attaching map for building $X_{l(i)+1}$ from $X_{l(i)}$. If $n_{r+i}-(r+i)= l(i)(p-1),$ which is in fact equivalent to condition (C) from \Cref{thm:total} in this case, consider the composition of maps of $p$-complete spectra as follows:
\begin{equation}\label[diagram]{eq:candidate}\begin{tikzcd}
\cp{r}{n} \ar[r]& A_i\ar[r]&\sphere^{2n_{r+i}} \ar[rr,"{\Sigma^{2(r+i)+1}\alpha_1^{(l(i))}}"] & & {\Sigma^{2(r+i)+1}X_{l(i)}}\ar[r]& \Sigma A_i \ar[r]& \Sigma \cp{r}{n}. \end{tikzcd} \end{equation} 

\begin{lemma}\label{lem:same-adams} \Cref{eq:candidate} gives a non-zero class after tensoring with $\eo$.
\end{lemma}
\begin{proof} 
To deduce the result, we first note that, for $n_{r+i}\equiv p-1 \pmod{p}$, the compositions \[\cp{r}{n} \to A_i \to \sphere^{2n_{r+i}}\] and 
\[\Sigma^{2(r+i)}X_{l(i)} \to A_i \to \cp{r}{n} \] are both $\eo$-split. So it suffices to note that the appropriate suspension of $\alpha_1^{(l(i))}$ is non-zero after tensoring with $\eo$, which follows from the proof of \Cref{prop:alphal}.
\end{proof}
\begin{ex}
If $r\equiv 1 \pmod{p}$ and $n=r+2p-2$, we are beyond the range where \Cref{geometric-easy} applies, but $\eo$-homology detects that there is a nontrivial, stably trivial $p$-torsion bundle $V$ pulled back from $S^{2n}$. If we take $N$ such that $r+3p-3>N>n=r+2p-2$, this bundle extends over the cellular inclusion $\cp{}{n} \to \cp{}{N}$, by simply noting that the $(p-1)$-Adams summand  is not altered until $N\geq r+3p-3$.
\end{ex}

For large primes $p$, the explicit construction of \Cref{eq:candidate} generally accounts for only some of the $p$-torsion bundles detected in $\eo$-theory, as we will demonstrate.

\begin{ex}\label{ex:7} Consider \Cref{thm:X} for $p= 7$. We find that $\Vect_{p^2}^0(\cp{}{p^2+(p-1)^2})$ surjects onto five summands isomorphic to $\Z/7$ in $\eo$-theory. If we take a generator for each summand in $\eo$-theory and lift it to a vector bundle representative, we find that only one generator can be represented by a simple extension of a vector bundle pulled back from a sphere (as in \Cref{eq:candidate}). In other words, take $V$ to be a non-trivial $7$-power-torsion bundle detected in one of the other summands. Consider a small-codimension cellular inclusion \[\iota_m \: \cp{}{p^2+(p-1)^2-m} \to \cp{}{p^2+(p-1)^2}\] for $0 \leq m \leq p-2$. We find that $\iota_m^*V$ is not pulled back from a sphere unless $\iota_m^*V$ is trivial. 
\end{ex}
\Cref{ex:7} contrasts with \Cref{thm:Y}: this example illustrates that, for larger corank, the $p$-torsion bundles detected by $\eo$-theory are not ones that we understand geometrically.

\section{%
    Computations of
    \texorpdfstring{$\eo_{-1}(\cp{r}{n}\otimes D \cp{r}{n})$}{EO−1(CP\_n\string^r ⊗ DCP\_n\string^r)}
    for
    \texorpdfstring{$n-r<\frac{2p^2-p-2}{2}$}{n-r < (2p²-p-2)/2}%
}
\label{combinatorics}

This section includes the proofs of technical results used in \Cref{lowerbound_computation}.
 We refer the reader to \Cref{subsec:alg_eo} for various definitions used throughout.
 The main content is the proof of \Cref{thm:total}.  Since this entire section is concerned with $\eo$-modules,
 the following definition will be useful. Throughout this section, $p$ is an odd prime. 

\begin{defn} With notation as in \Cref{defn:alg_eo_stuff}, let $X_l^e:= \eo\otimes X_l$. 
\end{defn}
Much of this section is most easily understood in terms of the combinatorics of the $(p-1)$ distinct $p$-complete Adams summands of the spectrum $\cp{}{\infty}$, or of the induced Adams summands of $\cp{r}{n}$. 

We recall that the spectrum $\cp{r}{n}$ admits a $p$-complete Adams decomposition $\cp{r}{n}\cong \oplus_{i=0}^{p-2}A_{i}$, where the spectrum $A_{i}$ consists of those cells in $\cp{r}{n}$ in degrees congruent to $2(r+i) \pmod{2p-2}$. Since each $A_k$ is $2p-2$-sparse, we can apply \cite[Theorem 5.12]{Chat} to deduce that, after tensoring with $\eo$, we have a further splitting of each $A_k$ as a sum of shifts of spectra of the form $X_l^e= \eo \otimes X_l$, where the spectra $X_l$ that occur in this decomposition are uniquely determined by the action of the Steenrod operation $P^1$ on the $\pmod{p}$-homology of $\cp{r}{n}$ (see \Cref{defn:alg_eo_stuff} for the definition of $X_l$). The group $\eo_{-1}\cp{r}{n} \otimes D\cp{r}{n}$ is then determined by where, and with what suspensions, the spectra $X^e_l$ occur in the $\eo$-module decomposition of each $\eo \otimes A_k$. Hence, it is essential for us to introduce notation to keep track of the degrees of ($\eo$-module) cells occurring in different Adams summands.

Recall the following definition.

\begin{defn} Given $k \in \mathbb Z$ and $d$ a positive integer, let $[k]_d=k \pmod d\in \{0,\ldots, d-1\}.$ For integers $n$ and $j$, let $n_j =n-[n-j]_{p-1}$, i.e., $n_j$ is the largest integer less than or equal to $n$ with the same remainder as $j$ $\pmod{p-1}$.
\end{defn}

\begin{defn}\label{def:I} Let $n,r$ be given with $n>r$. Let $p$ be prime. 
\begin{enumerate} 

\item Let 
\[l(n,r,i)=\op{min}\left\{ [r+i]_{p}+1, \frac{n_{r+i}-(r+i)}{p-1}+1\right\} .\]

When $n$ and $r$ are fixed, we will suppress dependence and write $l(n,r,i)=\l$.

\item Let $\mathcal I_{n,r}$  be the subset of $i \in \{0,\ldots, p-2\}$ such that 
\[ r+i + ([r+i]_p)(p-1) < n_{r+i} .\]
\item Let 
$y(n,r,i)=2n_{r+i}-2(p-[n_{r+i}]_p-1)(p-1)$. When $n$ and $r$ are fixed, we will suppress dependence and write $y(n,r,i)=\t$.
\end{enumerate}
\end{defn}
To understand the definitions above in terms of the combinatorics of Adams summands, note the following:
\begin{itemize}
\item  The spectrum $\Sigma^{2r+2i}X_{\l}$ is the term $B_0$ in the cofiltration of the $(r+i)$-Adams summand provided by \Cref{eq:filtfilt}; in other words, $\Sigma^{2r+2i}X_{\l}$ is the lowest-degree $\eo$-split summand in $\eo \otimes A_{i}$.
\item The set $\mathcal I_{n,r}$ indexes the Adams summands of $\cp{r}{n}$ for which the cofiltration given by \Cref{eq:filtfilt}
 has at least two terms; in other words, those Adams summands $A_{i}$ where $\eo \otimes A_{i}$ is a decomposable $\eo$-module. For $n-r$ large, $\mathcal I_{n,r}$ will include all $i$ between $0$ and $p-2$, inclusive. In particular, if $i\in \mathcal I_{n,r}$, then $l(n,r,i)=[r+i]_p+1$.
\item  In the case that $i \in \mathcal I_{n,r}$,  the top cofiltration quotient $B_t$ in the cofiltration of the $(r+i)$-Adams summand provided by \Cref{eq:filtfilt} 
is $\Sigma^{\t} X_{p-[n_{r+i}]_p}$. In other words, the top $\eo$-split summand of $\eo \otimes A_{i}$ is $ \Sigma^{\t} \eo \otimes X_{p-[n_{r+i}]_p}$.
\end{itemize}
We prove these observations next.
\begin{prop}\label{prop:cpnreo} Given integers $n>r \geq 0$, $\eo \otimes \cp{r}{n}$ has a summand

\[\bigoplus_{i=0}^{p-2} \Sigma^{\b} X^e_{\l} \oplus \bigoplus_{i\in \mathcal I_{n,r}} \Sigma^{\t } X^e_{\lp}.\]
The complementary summand splits as a sum of even suspensions of $X^e_p$. 
\end{prop}
\begin{proof} As noted in \Cref{ex:cpnr}, $\cp{r}{n}$ has algebraic $\eo$-theory. By \Cref{lem:hood_wants_many_lemmas},
 the decomposition of $\eo\otimes\cp{r}{n}$ 
into sums of even suspensions of $X_l$'s is determined by the 
$P^1$-module structure on the mod $p$-cohomology of 
$\cp{r}{n}$. The latter is, in turn, determined by the 
$P^1$-module structure on the mod $p$-cohomology of 
$\cp{}{\infty}$.

The degree $2k$ generator in $H^*(\cp{}{\infty};\Z/p)$ is $t^k$, where $t$
generates $H^2(\cp{}{\infty};\Z/p)$. Note that for $1\leq j  \leq [k]_p$,
\[(P^1)^j(t^k)=\left( \prod_{l=0}^{j-1}\left(k+(p-1)l\right)\right)t^{k+j(p-1)}\]
is nonzero and that 
\[(P^1)^{[k]_p+1}(t^k)=0.\]
For compatibility with \Cref{subsec:alg_eo}, we dualize and work with comodules over homology.
 Recall from \Cref{defn:alg_eo_stuff} that the quotient coalgebra $\F_p[\xi_1]/(\xi_1^p)$ of the dual Steenrod algebra 
has $p$ indecomposable graded comodules and
the comodule of length $i$ is denoted by $W_i$. The above computation shows that, as $\xi_1$-comodules,

\begin{equation}\label[empty]{eq:Adamsdecomp1}H_*(\cp{r}{\infty}) \cong \bigoplus_{i=0}^{p-2} 
\left ( \Sigma^{2r+2i}W_{[r+i]_p+1} \oplus \bigoplus_{j \geq 0} \Sigma^{2r+2i+ 2\left([r+i]_p +1 + jp\right)(p-1)} 
W_p \right ).\end{equation}
We will refer to the $i$-indexed $\xi_1$-comodule summands 
in \Cref{eq:Adamsdecomp1} above as Adams summands since they correspond to the Adams decomposition of $\cp{}{\infty}$ \cite{McG}.

Given a graded $\xi_1$-comodule $W$, let $(W)_k$ be the quotient of $W$ by elements of degree greater than or equal to $2k$. Then:
\begin{align*} H_*(\cp{r}{n}) &\cong H_*(\cp{r}{\infty})_{n+1}  \\
&\cong   \bigoplus_{i=0}^{p-2}
 \left (\Sigma^{2r+2i}\left(W_{[r+i]_p+1}\right)_{n+1} \oplus \bigoplus_{j \geq 0} 
\left(\Sigma^{2r+2i+ 2\left([r+i]_p +1 + jp\right)(p-1)} W_p\right)_{n+1} \right ). \end{align*}
Within each Adams summand, there is at most one summand of
\begin{equation}\label[empty]{eq:blehbleh}\left(\bigoplus_{j \geq 0} \Sigma^{2r+2i+ 2\left([r+i]_p +1 + jp\right)(p-1)} W_p\right)_{n+1}\end{equation}
that is neither zero nor a suspension of $W_p$.
 Since $H^*(\cp{r}{n})$ has a nonzero element in each even degree, 
the highest degree element in the Adams summand, with all elements in degrees congruent to 
$(r+i) \pmod{p-1}$ is the element in degree $n_{r+i}$.
In order for the module \Cref{eq:blehbleh} to be nonzero,
 the generator in degree $n_{r+i}$ must not lie in the bottom summand 
$\Sigma^{2r+2i}W_{[r+i]_p+1}$. Since $\Sigma^{2r+2i}W_{[r+i]_p+1}$ 
is concentrated in degrees $2(r+i)$ through 
$2\left(r+i+[r+i]_p(p-1)\right)$, this is equivalent to the condition that 
$n_{r+i}> r+i+[r+i]_p(p-1)$, i.e., that $i \in \mathcal I_{n,r}$.

In the case $i\in \mathcal I_{n,r}$, the interesting summand in the module \Cref{eq:blehbleh}
has top degree element in degree $n_{r+i}$ 
and contributes a summand $\Sigma^{ 2n_{r+i}-2(p-1)(p-[n_{r+i}]_p-1) }W_{p-[n_{r+i}]_p}$ to the top of the $(r+i)$-Adams summand.

If $i\notin \mathcal I_{n,r}$, note that 
\[\left(\Sigma^{2r+2i}W_{[r+i]_p+1}\right)_{n+1}\] 
has $\frac{n_{r+i}-(r+i)}{p-1}+1$ cells and begins in degree $2r+2i$, so that
\[\left(\Sigma^{2r+2i}W_{[r+i]_p+1}\right)_{n+1}\cong  \Sigma^{2r+2i}W_{\l}.\] 


Thus, we obtain the following:
\begin{align*} H_*(\cp{r}{n}) &\cong \left(H_*(\cp{r}{\infty})\right)_{n+1} \\
&\cong   \bigoplus_{i=0}^{p-2} \left( \Sigma^{2r+2i}W_{[r+i]_p+1} \oplus \bigoplus_{j \geq 0} \Sigma^{2r+2i+ 2\left([r+i]_p +1 + jp\right)(p-1)} W_p \right)_{n+1}\\
&\cong \left(\bigoplus_{i=0}^{p-2} \Sigma^{2(r+i)}W_{\l} \oplus \bigoplus_{i \in \mathcal I_{n,r}} \Sigma^{2n_{r+i} - 2(p-[n_{r+i}]_p-1)(p-1)}
W_{p-[n_{r+i}]_p} \right) \oplus T,
\end{align*}
where $T$ is a sum of shifts of $W_p$.
 \end{proof}

Our next ingredient is the $\eo$-module structure of $X_{l'}^e \otimes X_l^e$, necessary to describe $\eo\otimes\cp{r}{n} \otimes D\cp{r}{n}$.
\begin{lemma}[{\cite[Corollary 5.23]{Chat}}]\label{XiXj_decomp} Let $1 \leq l \leq l' \leq p$. 
Let $t=\begin{cases} l &\text{ if } l+l' \leq p \\
p-l' &\text{ if } l+l'>p.
\end{cases}.$
Then

\[X_l^e \otimes X_{l'}^e \simeq Z \oplus \bigoplus_{j=1}^t \Sigma^{2(l-j)(p-1)}X_{l'-l+2j-1}^e,\]

where $Z$ is a sum of even suspensions of $X_p^e$.
\end{lemma}
We also remind the reader that the $\eo$-homology of suspensions of $X_l$'s was computed in \Cref{prop:alphal}, which we restate below:
\begin{prop}\label{prop:alphal2}     
For $p$ an odd prime and $-\frac{(2p^2-p-2)}{2}< s < \frac{2p^2-p-2}{2}-(p-1)(l-1)$, $1\leq l \leq
p$,
\[\eo_{-1}\Sigma^{2s}X_l \cong \begin{cases} \Z/p \text{ if $l \neq p$ and $s=-l(p-1)$}\\ 0 \text{ otherwise.}\end{cases}\]
\end{prop}
Together, \Cref{prop:alphal2} and \Cref{XiXj_decomp} imply:
\begin{lemma}\label{lem:eo1zero} Let $l,l' \in \{1,\ldots, p\}$ with $l \leq l'$. Let $s$ be such that 
\[-\frac{(2p^2-p-2)}{2}< s < \frac{2p^2-p-2}{2}-(p-1)(l'+t-2),\] where $t$ is from \Cref{XiXj_decomp}. If either:
\begin{itemize} 
\item $s \not\equiv 0 \pmod{p-1}$, or
\item $l'=p$,
\end{itemize}
then $\eo_{-1}(\Sigma^{2s}X_l \otimes X_{l'})=0$.
\end{lemma}

Recall that \Cref{thm:total} is the main goal of this section. We state this theorem again for convenience, and then give its proof: 
\begin{thm}\label{thm:total2} Let $p$ be an odd prime. Let $n-r <\frac{2p^2-p-2}{2}.$ 
Then $\eo_{-1}(\cp{r}{n} \otimes D\cp{r}{n})\cong \left(\Z/p\right)^{\oplus j(n,r)}$ 
where $j(n,r)$ is the number of $i$ in $\{0,\ldots, p-2\}$ such that: 
\begin{itemize}
\item[{\bf A}.]  $ r+i + ([r+i]_p)(p-1) < n_{r+i} ,$ i.e., $i \in \mathcal I_{n,r}$;
\item[{\bf B}.]  $l(i)\leq [n_{r+i}]_p$; and
\item[{\bf C}.] $n_{r+i}-(r+i)+(p-1)\left([n_{r+i}]_p-[r+i]_p\right) = p(p-1).$
\end{itemize}
\end{thm}

\begin{proof}
We will write $Z \simeq_{X_p} Z'$ to indicate that two $\eo$-modules $Z$ and $Z'$
 become isomorphic after summing with even suspensions of $X_p^e$. 

By \Cref{prop:cpnreo}, there is an isomorphism of $\eo$-modules
\[\eo \otimes \cp{r}{n} \simeq_{X_p} \bigoplus_{\substack{ i=0 \\ \l\neq p}}^{p-2} 
\Sigma^{\b} X^e_{\l} \oplus \bigoplus_{\substack{i\in \mathcal I_{n,r}\\ [n_{r+i}]_p\neq 0}} \Sigma^{\t } X^e_{\lp}.\]
Since $DX_i \simeq \Sigma^{-2(i-1)(p-1)}X_i$,
\begin{align*}
\eo \otimes D\cp{r}{n}
& \simeq_{X_p}\,\,\eo \otimes \left(  \bigoplus_{\substack{{i=0}\\ \l \neq p}}^{p-2} \Sigma^{-\b} 
DX_{\l} \oplus \bigoplus_{\substack{ i\in \mathcal I_{n,r} \\ [n_{r+i}]_p \neq 0}} \Sigma^{-\t } DX_{\lp}\right) \\ 
& \simeq_{X_p}  \,\, \bigoplus_{\substack{{i=0}\\  \l \neq p}}^{p-2} \Sigma^{-\b-2(\l-1)(p-1)}X^e_{\l} \oplus \\ 
& \,\,\,\,\,\,\,\,\,\,\,\, \bigoplus_{\substack{i\in \mathcal I_{n,r}\\ [n_{r+i}]_p\neq 0}} \Sigma^{-\t-2(\lp-1)(p-1) }X^e_{\lp}\\ 
& \simeq_{X_p}\,\,  \bigoplus_{\substack{i=0\\  \l \neq p} }^{p-2} \Sigma^{-2\left(r+i+(\l-1)(p-1)\right)}
X^e_{\l} \oplus \bigoplus_{\substack{i\in \mathcal I_{n,r}\\ [n_{r+i}]_p\neq 0}} \Sigma^{-2n_{r+i} }X^e_{\lp}.
\end{align*}

To compute $\eo_{-1}\cp{r}{n} \otimes D \cp{r}{n}$, note that by \Cref{prop:alphal} and \Cref{lem:eo1zero} 

\begin{align}\label{eq:big}\eo\otimes\cp{r}{n}\!\otimes \! \! D\cp{r}{n}  &\simeq 
\left(\bigoplus_{\substack{i=0\\ \l \neq p}}^{p-2} 
\Sigma^{-2(\l-1)(p-1)} \left(X^e_{l(i)}\right)^{\otimes 2}\right) \, \oplus \, 
 \left(\bigoplus_{\substack{i \in \mathcal I_{n,r},\,
\l \neq p \\ [n_{r+i}]_p \neq 0 }} \!\!\!\Sigma^{2(r+i-n_{r+i} )} X_{\l}^e
\otimes X^e_{\lp} \right)\\
& \,\oplus \, \nonumber \left( \bigoplus_{\substack{ i \in \mathcal I_{n,r}\\ [n_{r+i}]_p \neq 0 } } 
\Sigma^{\t-2n_{r+i}} \left(X^e_{p-[n_{r+i}]_p}\right)^{\otimes 2} \right)\,
\\ & \oplus \, \nonumber
 \left( \bigoplus_{\substack{i \in \mathcal I_{n,r}, \, \l \neq p \\ [n_{r+i}]_p \neq 0 } } \Sigma^{\t-2(r+i+(\l-1)(p-1))} X_{\l}^e\otimes X^e_{\lp} \right)  \bigoplus \, C,
\end{align}
where $C$ is an $\eo$-module with $\pi_{-1}C=0$.

We now consider which summands of the previous formula can contribute to $\pi_{-1}$. 
\begin{itemize}
\item {\bf The fourth summand of \Cref{eq:big}}. Essentially for connectivity reasons, there are no contributions from the fourth summand. To work this out precisely, we consider two cases depending on whether $l\leq \lp$ or $l>\lp$.
\begin{itemize}
\item[{\bf Case 1}.] $\l \leq \lp$.  By \Cref{XiXj_decomp},

\[ X_{\l}^e \otimes X^e_{\lp} \simeq_{X_p} \bigoplus_{j=1}^t \Sigma^{2(\l-j)(p-1)} X^e_{\lp-\l+2j-1},\]
where 
\[
t = 
\begin{cases} 
    \l & \text{ if  } \l + \lp \leq p\\
    [n_{r+i}]_p& \text{ if  } \l + \lp >p.
\end{cases}
\]

Applying $\Sigma^{\t-2(r+i+(\l-1)(p-1))}$, we are computing a sum of terms of the form

\[\pi_{-1}\Sigma^{\t-2(r+i+(\l-1)(p-1)-(\l-j)(p-1))} X^e_{\lp-\l+2j-1},\]
for $1\leq j \leq t$.
The suspension in the $j$-th term is 
\[ R(i,j)= \t-2(r+i+(\l-1)(p-1)-(\l-j)(p-1) ).\]
\begin{claim} With $R(i,j)$ as above, $\left(  2p^2-p-2-2(p-1)(L-1)\right) >R(i,j)>0,$ where $L=\lp-\l+2j-1$. \end{claim}
Given the claim, $\pi_{-1}\Sigma^{R(i,j)} X^e_{\lp-\l+2j-1}=0$ by \Cref{prop:alphal2}. Therefore all groups that can contribute are actually zero.

To prove the claim, we use that $\t=2n_{r+i}-2(p-[n_{r+i}]_p-1)(p-1)$ to rewrite:
\begin{align*}
\frac{R(i,j)}{2(p-1)}=\frac{ n_{r+i}-(r+i)}{p-1} -(p-[n_{r+i}]_p+(j-2)).
\end{align*}
Note that:
\begin{itemize}
\item  $\frac{n_{r+i}-(r+i)}{p-1}$ is one less than the number of cells in the $(r+i)$-Adams summand of $\cp{r}{n}$; 
\item $p-[n_{r+i}]_p$ is the number of cells in the top $\eo$-split quotient of the $(r+i)$-Adams summand of $\cp{r}{n}$;
\item $j-2$ is less than $l(i)-1=[r+i]_p$, which is one less than the number of cells in the bottom $\eo$-split subobject of the $(r+i)$-Adams summand of $\cp{r}{n}$. 
\end{itemize}

In particular, since $i\in \mathcal I_{n,r}$, the top $\eo$-split quotient and the bottom $\eo$-split subobject in the $(r+i)$-Adams summand consist of disjoint cells, and $p-[n_{r+i}]_p+(j-2)< \frac{n_{r+i}-(r+i)}{p-1}$. 

The upper bound on $R(i,j)$ simply follows from noting that $\Sigma^{R(i,j)}X_{\lp-\l+2j-1}$ is an $\eo$-split summand in $\eo \otimes \cp{r}{n} \otimes D \cp{r}{n}$, whose top cell is in dimension less than $(2p^2-p-2)$.

 \item[{\bf Case 2}.] $\lp<\l$.  By \Cref{XiXj_decomp},

\[X_{\l}^e \otimes X^e_{\lp} \simeq_{X_p} \bigoplus_{j=1}^t \Sigma^{2(\lp-j)(p-1)} X^e_{\l-(\lp)+2j-1},\]
where \[t = \begin{cases} \lp& \text{ if  } \l + \lp \leq p\\
p-\l& \text{ if  } \l + \lp >p. \end{cases}\]

Applying $\Sigma^{\t-2(r+i+(\l-1)(p-1))}$, we are computing a sum of terms of the form

\[\pi_{-1}\Sigma^{\t-2(r+i+(\l-1)(p-1)-  (\lp-j)(p-1)   )} X^e_{\l-(\lp)+2j-1},\]
for $1\leq j \leq t$.
The suspension in the $j$-th term is 
\[ S(i,j)= \t-2(r+i+(\l-1)(p-1)-(\lp-j)(p-1) ).\]
\begin{claim} With $S(i,j)$ as above, $\left(  2p^2-p-2-2(p-1)(L-1)\right) >S(i,j)>0,$ where $L=\l-(\lp)+2j-1$.\end{claim}
Given the claim, $\pi_{-1}\Sigma^{S(i,j)} X^e_{\l-(\lp)+2j-1}=0$ by \Cref{prop:alphal2}. Therefore all groups that can contribute are actually zero.

To prove the claim, we use that $\t=2n_{r+i}-2(p-[n_{r+i}]_p-1)(p-1)$ to rewrite:
\begin{align*}
\frac{S(i,j)}{2(p-1)}=\frac{ n_{r+i}-(r+i)}{p-1} -((\l-1)+(j-1)).
\end{align*}
Note that:
\begin{itemize}
\item  $\frac{n_{r+i}-(r+i)}{p-1}$ is one less than the number of cells in the $(r+i)$-Adams summand of $\cp{r}{n}$; 
\item $\l-1$ is one fewer than the number of cells in the bottom $\eo$-split subobject of the $(r+i)$-Adams summand of $\cp{r}{n}$.  
\item $j-1$ is at most $t-1$, and $t$ is at most $p-[n_{r+i}]_p$, which is the number of cells
in the top $\eo$-split quotient of the $(r+i)$-Adams summand of $\cp{r}{n}$.
\end{itemize}
In particular, since $i\in \mathcal I_{n,r}$, the top $\eo$-split quotient and the bottom $\eo$-split subobject in the $(r+i)$-Adams summand consist of disjoint cells in the $(r+i)$-Adams summand, and $(\l-1)+(j-1)< \frac{n_{r+i}-(r+i)}{p-1}$. 

The upper bound on $S(i,j)$ simply follows from noting that $\Sigma^{S(i,j)}X_{\l-(\lp)+2j-1}$ is an $\eo$-split summand in $\eo \otimes \cp{r}{n} \otimes D \cp{r}{n}$, whose top cell is in dimension less than $2p^2-p-2$.
\end{itemize}

 \item {\bf The first and third summands of \Cref{eq:big}.} For both we will use the formula
\begin{equation}\label[empty]{eq:Xmm}X^e_{m}\otimes X^e_{m} \simeq_{X_p}  
\bigoplus_{j=1}^t \Sigma^{2(m-j)(p-1)}X^e_{2j-1},\end{equation}
where $t=\op{min}\{m,p-m\}$, which follows from taking $l=l'=m$ in \Cref{XiXj_decomp}.
Combining \Cref{eq:Xmm} with \Cref{prop:alphal2}, we see that
\[\pi_{-1} \Sigma^{2R}(X_m^e \otimes X_m^e) \neq 0\]
if and only if there exists $1 \leq j \leq \op{min}\{m,p-m\}$ such that $2j-1 \neq p$ and
\begin{equation}\label[empty]{eq:condition00}-(2j-1)(p-1)=R+(m-j)(p-1).\end{equation}
We consider the possible values for $R$ and $j$ 
that occur in the first and third summands of \Cref{eq:big}.
\begin{itemize}
\item {\bf First summand.} In this case, 
$R= -(\l-1)(p-1)$ and $m=\l$. Substituting in 
\Cref{eq:condition00}, we must solve for $j$ such that 
\[-(2j-1)(p-1)=-(m-1)(p-1)+(m-j)(p-1).\] 
Dividing by $p-1$, $-2j+1=-j+1.$ Since $1 \leq j \leq \op{min}\{m, p-m\}$, there is no solution.
\item {\bf Third summand.} In this case, \[2R= \t-2n_{r+i}=-2(p-[n_{r+i}]_p-1)(p-1),\]
and \[R=([n_{r+i}]_p-p+1)(p-1).\] Moreover, $m=p-[n_{r+i}]_p$, so we get that
\[R=(-m+1)(p-1).\] 
Substituting in \Cref{eq:condition00}, 
we must solve for $j$ such that 
\[-(2j-1)(p-1)=(-m+1)(p-1)+(m-j)(p-1).\]
Dividing by $(p-1)$, we again find $-2j+1=-j+1$, and
there is no solution.
\end{itemize} 

\item {\bf The second summand of \Cref{eq:big}}. The argument above shows that only terms from the second summand in \Cref{eq:big} can contribute, which actually provides clarity on the possible maps $\cp{r}{n} \to \Sigma \cp{r}{n}$ which are non-zero after tensoring with $\eo$:
\begin{cor}\label{lem:same-adams-0} The terms in $\eo_{-1}$ of the second summand of \Cref{eq:big} correspond to maps from the top 
$\eo$-split subcomplex within a given Adams summand of $\cp{r}{n}$ to the bottom 
$\eo$-split subcomplex within the same Adams summand. In particular, as maps from $\eo \otimes \cp{r}{n} $ to $\eo \otimes \Sigma \cp{r}{n}$, they are the Hurewicz image of maps of the form
\begin{equation}\label{eq:adams-factor-0}
\cp{r}{n} \to A_i\to \Sigma A_i \to \Sigma \cp{r}{n}, \end{equation}
where $A_i$ is the Adams summand in $\cp{r}{n}$ consisting of cells in degrees $2k$ for which $k \equiv r+i \pmod{p-1}$. 
 \end{cor}
Note that we have shown that only the second summand of \Cref{eq:big} can {\em possibly} contribute to the $\eo_{-1}$-cohomology, but to determine which terms actually contribute requires more work. We consider various cases to apply \Cref{XiXj_decomp}.
\begin{itemize}
\item[]
\begin{enumerate}
\item[{\bf Case 1}.] $\l \leq \lp$.  By \Cref{XiXj_decomp},

\[ X_{\l}^e \otimes X^e_{\lp} \simeq_{X_p} \bigoplus_{j=1}^t \Sigma^{2(\l-j)(p-1)} X^e_{\lp-\l+2j-1},\]
where 
\[
t = 
\begin{cases} 
    \l & \text{ if  } \l + \lp \leq p\\
    [n_{r+i}]_p& \text{ if  } \l + \lp >p.
\end{cases}
\]
Therefore,
\begin{equation}
    \label[empty]{eq:tb1}   
    \Sigma^{2(r+i-n_{r+i} )} X_{\l}^e\otimes X^e_{\lp} \simeq_{X_p} \bigoplus_{j=1}^t \Sigma^{2(r+i)-2n_{r+i} +2(\l -j)(p-1)}X^e_{\lp-\l+2j-1}.
\end{equation}
By \Cref{prop:alphal2}, we see that \Cref{eq:tb1} contributes a $\Z/p$ to $\eo_{-1}\cp{r}{n} \otimes D\cp{r}{n}$ for each $1 \leq j \leq t$ such that
\[
    r+i-n_{r+i} +( \l -j)(p-1) = -\left(\lp - \l +2j-1\right)(p-1). 
\]
This can be rewritten as:
\begin{equation}\label[empty]{eq:tb2} 
n_{r+i}-(r+i)= \left(\lp +j-1\right)(p-1). 
\end{equation}

Note that $r+i\equiv n_{r+i}\pmod{p-1}$. Rearranging and dividing by $(p-1)$, we get:

\begin{equation}\label[empty]{eq:tb2-1} 
[n_{r+i}]_p+\frac{n_{r+i}-(r+i)}{p-1}=p+j-1. 
\end{equation}
Note that $[n_{r+i}]_p+\frac{n_{r+i}-(r+i)}{p-1}\equiv [r+i]_p \pmod{p}$, which implies that $j\equiv [r+i]_p+1\pmod{p}$. Note that $i\in \mathcal I_{n,r}$ implies that $l(i)=[r+i]_p+1$. Thus the only candidate solution to to \Cref{eq:tb2} with $1\leq j \leq p-1$ is $j=[r+i]_p+1$.

However, we also need to satisfy the possibly stronger condition, $1\leq j \leq t$. If $l(i)\leq [n_{r+i}]_p$, then $t=l(i)$ and $j=[r+i]_p+1=l(i)$ is a valid choice. If $l(i)>[n_{r+i}]_p$, then there is no solution.

All in all, we find that those $i$ that contribute a summand of $\Z/p$ to the computation are exactly those (A) $i\in \mathcal I_{n,r}$ so that (B) $[r+i]_p< [n_{r+i}]_p$, and (C):
\[
n_{r+i}-(r+i)= \left(\lp +[r+i]_p\right)(p-1),
\]
which we obtain by substituting $j=[r+i]_p+1$ in \Cref{eq:tb2}.

 \item[{\bf Case 2}.] $\lp<\l$.  By \Cref{XiXj_decomp},

\begin{equation}\label[empty]{eq:xixjcase2} X_{\l}^e \otimes X^e_{\lp} \simeq_{X_p} \bigoplus_{j=1}^t \Sigma^{2(\lp-j)(p-1)} X^e_{\l-(\lp)+2j-1},\end{equation}
where \[t = \begin{cases} \lp& \text{ if  } \l + \lp \leq p\\
p-\l& \text{ if  } \l + \lp >p. \end{cases}\]
Applying $\Sigma^{2(r+i-n_{r+i})}$ to both sides of \Cref{eq:xixjcase2} and appealing to \Cref{prop:alphal2}, we see that summands in \Cref{eq:xixjcase2} contribute to $\eo_{-1}\cp{r}{n} \otimes D\cp{r}{n}$ if and only if 
\[ (\lp-j)(p-1)+(r+i)-n_{r+i}=-\left( \l-(\lp)+2j-1\right)(p-1),\]
which rearranges to:
\[n_{r+i}-(r+i)=(j+ \l-1)(p-1).\]
Using that $i\in \mathcal I_{n,r}$ so that $\l=[r+i]_p+1$, this reduces to:
\begin{equation}\label[empty]{eq:tb7} n_{r+i}-(r+i)=(j+ [r+i]_p)(p-1). \end{equation} 
Reducing both sides of \Cref{eq:tb7} modulo $p$ gives that $j\equiv -[n_{r+i}]_p\pmod{p}.$ Since $1 \leq j \leq p-1$, the only possible solution is $j=p-[n_{r+i}]_p$. However, $j$ must also satisfy the stronger condition $1 \leq j \leq t$. 

If $l(i)\leq [n_{r+i}]_p$, then $t=p-[n_{r+i}]_p$, so this solution is possible. If $l(i)>[n_{r+i}]_p,$ then $t=p-l(i)$ and this choice of $j$ is not possible.

All in all, we find that those $i$ that contribute a summand of $\Z/p$ to the computation are exactly those (A) $i\in \mathcal I_{n,r}$ so that (B) $[r+i]_p< [n_{r+i}]_p$, and (C):
\[
n_{r+i}-(r+i)= \left(\lp +[r+i]_p\right)(p-1),
\]
which we obtain by substituting $j=p-[n_{r+i}]_p$ in \Cref{eq:tb7}.

\end{enumerate}
\end{itemize}
\end{itemize}
\end{proof}


\bibliographystyle{abbrv}
\bibliography{eo_cplx_VBS}

@article{AdamsWalker,
  author={J. F. Adams and G. Walker},
  title={On complex {S}tiefel manifolds},
  journal={Proc. Camb. Phil. Soc.},
  volume={61},
  pages={81--103},
  year={1965},
}

@article{AR,
  author={M. Atiyah and E. Rees},
  title={Vector bundles on projective 3-space},
  journal={Inventiones Math.},
  volume={35},
  pages={131--153},
  year={1976},
}

@article{Arone02,
  author={G. Arone},
  title={The {W}eiss derivatives of {BO(-) and BU(-)}},
  journal={Topology},
  volume={41},
  number={3},
  pages={451--481},
  year={2002},
}

@article{Bousfield,
  author={A. Bousfield},
  title={The localization of spectra with respect to homology},
  journal={Topology},
  volume={18},
  number={4},
  pages={257--281},
  year={1979},
}

@misc{Chat,
  author={H. Chatham},
  title={An Orientation Map for Height $p-1$ Real {$E$} Theory},
  year={2019},
  howpublished={arxiv:1908.11496},
}

@article{Devinatz,
  author={E. Devinatz},
  title={A {L}yndon--{H}ochschild--{S}erre spectral sequence for certain homotopy fixed point spectra},
  journal={Trans. Amer. Math. Soc.},
  volume={357},
  number={1},
  pages={129--150},
  year={2005},
}

@article{Hu,
  author={Y. Hu},
  title={Metastable complex vector bundles over complex projective spaces},
  journal={Trans. Amer. Math. Soc.},
  volume={376},
  number={11},
  pages={7783--7814},
  year={2023},
}

@article{McG,
  author={C.A. McGibbon},
  title={Stable properties of rank 1 loop structures},
  journal={Topology},
  volume={20},
  pages={109-118},
  year={1981},
  publisher={Pergamon Press},
}

@inproceedings{Mimura_HBAT,
  author={M. Mimura},
  editor={I. M. James},
  title={Homotopy theory of {L}ie groups},
  booktitle={Handbook of Algebraic Topology},
  volume={58},
  pages={951--991},
  year={1995},
  publisher={Elsevier},
}

@article{MO,
  author={Morgan Opie},
  title={A classification of complex rank 3 vector bundles on {$\mathbb CP^5$}},
  journal={Advances in Mathematics},
  volume={455},
  pages={109878},
  year={2024},
  doi={https://doi.org/10.1016/j.aim.2024.109878},
  url={https://www.sciencedirect.com/science/article/pii/S0001870824003931},
  issn={0001-8708},
}

@article{mosher69,
  author={Mosher, Robert E},
  title={Some homotopy of stunted complex projective space},
  journal={Illinois Journal of Mathematics},
  volume={13},
  number={1},
  pages={192--197},
  year={1969},
  publisher={Duke University Press},
}

@article{Nave,
  author={L. S. Nave},
  title={The {S}mith--{T}oda complex $V((p+1)/2)$ does not exist},
  journal={Ann. of Math.},
  volume={171},
  number={1},
  pages={419--509},
  year={2010},
}

@article{opie24enum,
  author={Opie, Morgan P},
  title={Enumerating complex rank {$ n $} vector bundles on {$\mathbb C P^{n+1}$}},
  journal={arXiv preprint arXiv:2410.23520},
  year={2024},
}

@book{Ravenel_green,
  author={D. C. Ravenel},
  title={Complex cobordism and the stable homotopy groups of spheres},
  year={2004},
  publisher={AMS Chelsea Publishing},
}

@article{Toda59,
  author={H. Toda},
  title={A topological proof of theorems of {B}ott and {B}orel--{H}irzebruch for homotopy groups of unitary groups},
  journal={Memoirs of the College of Sciences, University of Kyoto},
  series={A},
  volume={XXXII},
  number={1},
  pages={103--119},
  year={1959},
}

@article{Weiss,
  author={M. Weiss},
  title={Orthogonal calculus},
  journal={Trans. Amer. Math. Soc.},
  volume={347},
  number={10},
  pages={3743--3796},
  year={1995},
}

\end{document}